\documentclass[sn-mathphys,Numbered]{sn-jnl}

\usepackage{graphicx}%
\usepackage{multirow}%
\usepackage{amsmath,amssymb,amsfonts}%
\usepackage{amsthm}%
\usepackage{mathrsfs}%
\usepackage[title]{appendix}%
\usepackage{xcolor}%
\usepackage{textcomp}%
\usepackage{manyfoot}%
\usepackage{booktabs}%
\usepackage{algorithm}%
\usepackage{algorithmicx}%
\usepackage{algpseudocode}%
\usepackage{listings}%
\usepackage{bbm}
\usepackage{setspace}
\usepackage{subcaption}
\usepackage{bm}

\theoremstyle{thmstyleone}%
\theoremstyle{thmstyletwo}%

\theoremstyle{thmstylethree}%

\begin{document}

\title[A Parallel-in-time Method Based on Preconditioner for Biot's Problem]{A Parallel-in-time Method Based on Preconditioner for Biot's Model}


\author[1]{\fnm{Zeyuan} \sur{Zhou}}\email{zeyuanzhou117@gmail.com}

\author[2]{\fnm{Huipeng} \sur{Gu}}\email{12131226@mail.sustech.edu.cn}

\author[3]{\fnm{Guoliang} \sur{Ju}}\email{jvguoliang@foxmail.com}

\author*[4]{\fnm{Wei} \sur{Xing}}\email{xingwei@nsccsz.cn}

\affil[1]{\orgdiv{Faculty of Aerospace Engineering}, \orgname{Technische Universiteit Delft}, \orgaddress{\city{Delft}, \postcode{2629 HS}, \country{Netherlands}}}

\affil[2]{\orgdiv{Department of Mathematics}, \orgname{Southern University of Science and Technology}, \orgaddress{\city{Shenzhen}, \postcode{518055},  \country{China}}}

\affil[3]{\orgdiv{Research and Development Department}, \orgname{Tenfeng}, \orgaddress{\city{Shenzhen}, \postcode{518055}, \country{China}}}

\affil*[4]{\orgdiv{High Performance Computing Department}, \orgname{National Supercomputing Center}, \orgaddress{ \city{Shenzhen}, \postcode{518055}, \country{China}}}


\abstract{We proposed a parallel-in-time method based on preconditioner for Biot's consolidation model in poroelasticity. In order to achieve a fast and stable convergence for the matrix system of the Biot's model, we design two preconditioners with approximations of the Schur complement. The parallel-in-time method employs an inverted time-stepping scheme that iterates to solve the preconditioned linear system in the outer loop and advances the time step in the inner loop. This allows us to parallelize the iterations with a theoretical parallel efficiency that approaches 1 as the number of time steps and spatial steps grows. We demonstrate the stability, accuracy, and linear speedup of our method on HPC platform through numerical experiments.}





\maketitle

\section{Introduction}\label{sec1}
The behavior of fluid-saturated porous materials undergoing deformation can be described by coupled poroelastic equations, which were first studied by Terzaghi \cite{Terzaghi}. Biot \cite{Biot1,Biot2} then extended the theory to a more general three-dimensional model. Biot’s model originated from geophysical applications, but also applied in various fields such as biomechanics, marine engineering, and petroleum engineering.

Poroelastic problem is commonly involved in numerical simulation of the multi-physics phenomenon, where the large-scale linear system resulted from the Biot's model is one of the most computation consuming part of the problem. For this reason, a lot of effort are made to designing efficiently parallel method for solving Biot's model. 
Temporal parallelization method is a promising strategy that gives a view to parallelize the computing in temporal dimension and becomes increasingly attractive in recent years. A famous kind of temporal parallelization method is known as parallel-in-time integration method, which uses two sets of temporal discretization grids and propagators - one fine and one coarse. The coarse propagator runs first sequentially, and then the fine propagators run in parallel to correct the coarse solution until the final solution converges to the sequential high-order time-stepping method such as the Runge-Kutta method. Several well-known methods include the Parareal method \cite{parareal}, the PFASST method \cite{PFASSt}, the MGRIT method \cite{MGRIT}, and the STMG method \cite{STMG}. Since the time stepping on the coarse grid is sequential in parallel-in-time integration methods, it can only parallelize the computation partially. Besides, another temporal parallelization strategy is the parallel-in-time fixed-stress splitting method, which splits the Biot's model into two steps and parallelize the fixed-stress iterations \cite{PFSparallel}.


In this paper, we propose a parallel-in-time method based on preconditioner for Biot's model. We consider a quasi-static Biot’s model as a time-dependent two-field system of displacement $\boldsymbol{u}$ and fluid pressure $p$ and discretized with P1-P1 finite element method, which leads to a typical saddle point problem \cite{gasparpeturb,gasperpeturb2}. 
Following the idea given by Murphy et al. \cite{precondition1} and Ipsen \cite{precondition2}, we design two preconditioners using an approximate Schur complementuse. These preconditioners ensure a stable convergence for the parallel-in-time method. Then, the parallel-in-time method based on preconditioner uses the inverted time-stepping method, which takes the time-stepping process as the inner loop while the iterations solving preconditioned linear system as the outer loop, and each thread only runs a very few iterations at each time step. This method requires no sequential initialization, and the theoretical parallel efficiency converges to $1$ as the number of time steps and spatial steps increase.

The remained part of the paper is organized as follows. In Section \ref{sec:modelanddiscre}, we introduce the quasi-static Biot's model and the stabilized P1-P1 finite element method used in this paper; in Section \ref{sec:pitbP}, the parallel-in-time method based on preconditioner and its implementation on HPC platform is proposed; in Section \ref{chapter_result}, several benchmarking numerical cases are used to study the consistency, stability, accuracy, and speedup ratio of the parallel method. Section \ref{sec:conclusion} gives the conclusion and recommendation for the future work based on the study results.

\section{Mathematical model and discretization}\label{sec:modelanddiscre}
In this section, we introduce the setting for the quasi-static biot's model and the discretization schemes in space and time, which are used in the parallel-in-time method based on preconditioner proposed in this paper. 
\subsection{Quasi-static biot's model}

We consider the quasi-static Biot's model in the two-field formulation on a domain $\Omega \subset \mathbbm{R}^2$ with boundary denoted as $\partial \Omega$. The process is described by the partial differential equations in space-time domain $\Omega \times(0,T]$ as follows:

\begin{align}
    & -\nabla \cdot \sigma(\bm{u})+ \alpha \nabla p =\boldsymbol{f}, \label{biotsys1}\\
    & \partial_t (\alpha \nabla \cdot {\bm{u}}) - K \Delta p  = g, \label{biotsys2}
\end{align}
where
\begin{align*}
	\sigma(\bm{u})=2\mu\varepsilon(\bm{u})+\lambda (\mbox{div}\bm{u}) \mathbf{I}, \ \ \ \varepsilon(\bm{u}) = \frac{1}{2} \left[ \nabla \bm{u} + (\nabla \bm{u})^T \right].
\end{align*}
Here, the displacement $\boldsymbol{u}$ and the fluid pressure $p$ are primary unknowns, $\alpha$ is the Biot-Wills constant, $K$ is the ratio of material permeability to fluid viscosity. $\mathbf{I}$ is the identity tensor, $\boldsymbol{f}$ is body force density and $g$ is the source term which describe the injection or extraction process of fluid. Lam{\'e} coefficients $\lambda$ and $\mu$ can be computed by the Young's modulus $E$ and the Poisson ratio $\nu$:
\begin{align*}
	\lambda = \frac{E\nu}{(1+\nu)(1-2\nu)},\ \ \ \mu = \frac{E}{2(1+\nu)}.
\end{align*}
Appropriate boundary and initial conditions should be provided to make the problem well-posed. For simplicity, we consider the following boundary conditions:
\begin{align}
    p =0, \quad {\rm on}\ \Gamma_{p,D} \times(0,T], 
    \\ 
    \sigma(\boldsymbol{u}) \boldsymbol{n}= \boldsymbol{0}, \quad \text{on} \ \Gamma_{p,N} \times(0,T], 
    \\
    \boldsymbol{u} = \boldsymbol{0}, \quad {\rm on}\ \Gamma_{\boldsymbol{u},D} \times(0,T], 
    \\
    (\nabla p)\cdot \boldsymbol{n} = 0, \quad \text{on} \ \Gamma_{\boldsymbol{u},N} \times(0,T],
\end{align}
where $\boldsymbol{n}$ is the unit outward normal to the boundary, $\Gamma_{p,D} \cup \Gamma_{p,N} = \Gamma_{\boldsymbol{u},D} \cup \Gamma_{\boldsymbol{u},N} = \partial \Omega$ with $\lvert \Gamma_{\bm{u},D} \rvert>0, \lvert \Gamma_{p,D} \rvert >0$. Next, we endow the problem with compatible the initial data:
\begin{align}
    p(0) = p^0, \quad \text{in} \ \Omega \times \{0\}, \label{intial1}\\
    \boldsymbol{u}(0) = \boldsymbol{u}^0, \quad \text{in} \ \Omega \times \{0\} \label{intial2},
\end{align}
The existence and uniqueness of solution of Biot's model \eqref{biotsys1}-\eqref{biotsys2} with initial condition \eqref{intial1}-\eqref{intial2} is studied in the works of Zenisek \cite{vzenivsek1984existence} and Showalter \cite{showalter2000diffusion}.

First, we introduce two Sobolev spaces as follows:
\begin{align}
    Q =& \left\{p \in H^1(\Omega),\ p =0 \ \  \text{on}\ \Gamma_{p,D}  \right\}, 
    \\
    \boldsymbol{V} =& \left\{{\boldsymbol{u} \in (H^1(\Omega))^2,\ \boldsymbol{u}=\boldsymbol{0} \ \ \text{on}\ \Gamma_{\boldsymbol{u},D } }
\right\}.
\end{align}
For ease of presentation, we define the following bilinear forms for $p,\psi \in Q$, and $\bm{u}, \bm{v} \in \boldsymbol{V}$:
\begin{align*}
    & a(\bm{u},\bm{v}) = 2\mu \int_\Omega  \varepsilon (\bm{u}) :  \varepsilon (\bm{v}),  \\
    & b(\bm{v},\phi) =  \alpha \int_\Omega \phi \  \nabla \cdot \bm{v}, \\
    & c(p,\psi) =  K \int_\Omega \nabla p \cdot \nabla \psi.
\end{align*}
Multiplying \eqref{biotsys1}-\eqref{biotsys2} by test functions, integrating by parts, and applying boundary conditions \eqref{intial1}-\eqref{intial2} yields the following variational formulation: for a given $t>0$, find $(\bm{u},p) \in \bm{V} \times Q$ such that
 \begin{align}
     a(\bm{u},\bm{v}) - b(\bm{v},p) = (\bm{f},\bm{v}), \quad  \forall \bm{v}\in \bm{V}, \label{eq:1} 
     \\
     b(\partial_t \bm{u},p) + c(p,q) = (g,q), \quad \forall q \in Q.  \label{eq:2}
 \end{align}

\subsection{Discretization scheme}\label{sec:stabfem}

It is necessary for a finite element pair of spaces $\boldsymbol{V}_h\times Q_h$ satisfying an inf-sup condition. We choose stabilized P1-P1 scheme, which is introduced in \cite{gasparpeturb} and analyzed in \cite{gasperpeturb2}, as the spatial discretizaition scheme. Denote the function spaces spanned by $C^0$ piecewise polynominal of 1st degree as $S^1_h\in H^1(\Omega)$, then in $\mathbbm{R}^2$, there are discretized spaces $\boldsymbol{V}_h = \boldsymbol{V}\cap (S^1_h\times S^1_h)$ and $Q_h = Q\cap S^1_h$. The semi-discrete Galerkin formulations of Eqs. \eqref{eq:1}  and \eqref{eq:2} can be derived as following:

Find $(\boldsymbol{u}_h(t),p_h(t))\in C^1([0,T];\boldsymbol{V}_h)\times C^1([0,T],Q_h)$ such that
\begin{equation}\label{eq:hformulat}
   \begin{aligned}
    a(\boldsymbol{u}_h,\boldsymbol{v}_h)-b( \boldsymbol{v}_h,p_h) &= (\boldsymbol{f}_h,\boldsymbol{v}_h), \ \ \forall \boldsymbol{v}_h\in\boldsymbol{V}_h, \ \ t\in (0,T],
    \\
  b(\partial_t {\boldsymbol{u}}_h,q_h)+c(p_h,q_h)&=(g_h,q_h), \ \ \forall q_h \in Q_h, \ \ t\in (0,T].
  \end{aligned} 
\end{equation}
 Stabilized P1-P1 scheme introduces a perturbation term $\beta \Delta ( \partial_t p )$ into primal system \eqref{eq:hformulat}, and the semi-discrete system can be written as:
\begin{equation}
   \begin{aligned}
    a(\boldsymbol{u}_h,\boldsymbol{v}_h)-b( \boldsymbol{v}_h, p_h) &= (\boldsymbol{f}_h,\boldsymbol{v}_h), \ \ \forall \boldsymbol{v}_h\in\boldsymbol{V}_h,
    \\
  b(\partial_t {\boldsymbol{u}}_h,q_h)+  \frac{\beta}{K}c(\partial_t {p}_h,q_h) +c(p_h,q_h)&=(g_h,q_h), \ \ \forall q_h \in Q_h.
  \end{aligned} 
\end{equation}
The factor $\beta$ can be calculated as:
\begin{equation}
    \beta = \frac{h}{4(\lambda +2\mu)},
\end{equation}
where $h$ is the mesh size, as for nonuniform mesh, generally the maxima is chosen. This strategy has been theoretically proven to be effective, for the details of proof, see in \cite{gasperpeturb2}.

Use backward-Euler time-stepping method for temporal discertization and set the time step length is $\tau$, then the fully discretized problem can be written as following: For $m\geq 1$ and $\left( 
  \boldsymbol{u}_h,p_h\right)\in \boldsymbol{V}_h\times Q_h$ there is
  \begin{equation}\label{eq:diseq1}
      \begin{aligned}
    a(\boldsymbol{u}^m_h,\boldsymbol{v}_h)-b(\boldsymbol{v}_h,p^m_h) = (\boldsymbol{f}^m_h,\boldsymbol{v}_h), \ \ \forall \boldsymbol{v}_h\in\boldsymbol{V}_h,
    \\
  b(\boldsymbol{u}^m_h,q_h)+(\tau  +\frac{\beta}{K})c(p^m_h,q_h)=b(\boldsymbol{u}^{m-1}_h,q_h)\\+\frac{\beta}{K} c(p^{m-1}_h,q_h)
  + \tau (g^m_h,q_h), \ \ \forall q_h  \in Q_h,
\end{aligned}
  \end{equation}
where $\tau$ is the time step size.

We denote the elasticity matrix as $A$, the gradient matrix as $B$ and the diffusive matrix as $C$. The relation of bilinear operators to matrices is following:
\begin{equation*} 
\begin{aligned}
    a(\boldsymbol{u}_h,\boldsymbol{v}_h) \rightarrow 
A, \\
-b(\boldsymbol{v}_h,p_h) \rightarrow B,\\
c(p_h,q_h) \rightarrow C,
\end{aligned}
\end{equation*}
then the matrix formulation of Eq. \eqref{eq:diseq1} can be derived as following:
\begin{equation}
\left[\begin{tabular}{cc}
     $A$& $ B$ \\
    $B^{\rm T}$ & $-(\tau  +\frac{\beta}{K} )C$ 
\end{tabular} \right]
\left[
\begin{tabular}{c}
     $\boldsymbol{U}^m$ \\
     $\boldsymbol{P}^m$ 
\end{tabular}
\right] = 
\left[\begin{tabular}{cc}
     $0$& $0$ \\
    $B^{\rm T}$ & $-\frac{\beta}{K} C$ 
\end{tabular} \right]
\left[
\begin{tabular}{c}
     $\boldsymbol{U}^{m-1}$ \\
     $\boldsymbol{P}^{m-1}$ 
\end{tabular}
\right] +
\left[
\begin{tabular}{c}
     $\boldsymbol{F}^m$ \\
     -$\tau\boldsymbol{G}^m$ 
\end{tabular}
\right].
\label{eq:fullysystem}
\end{equation}

 For convenience, In the following sections, Eq. \eqref{eq:fullysystem} is simplified as
\begin{equation}\label{eq:axb}
    \mathbbm{A}X =b,
\end{equation}
where
$$
\begin{aligned}
   \mathbbm{A} =  \left[\begin{tabular}{cc}
     $A$& $B$ \\
    $B^T$ & $-(\tau +\frac{\beta}{K} )C$ 
\end{tabular} \right], \quad 
X =\left[
\begin{tabular}{c}
     $\boldsymbol{U}^m$ \\
     $\boldsymbol{P}^m$ 
\end{tabular}
\right], \\
b = \left[\begin{tabular}{cc}
     $0$& $0$ \\
    $B^T$ & $-\frac{\beta}{K} C$ 
\end{tabular} \right]
\left[
\begin{tabular}{c}
     $\boldsymbol{U}^{m-1}$ \\
     $\boldsymbol{P}^{m-1}$ 
\end{tabular}
\right] +
\left[
\begin{tabular}{c}
     $\boldsymbol{F}^m$ \\
     $-\tau\boldsymbol{G}^m$ 
\end{tabular}
\right].
\end{aligned}
$$
For convenience in following sections, we also denote 
$$
f^m = \left[
\begin{tabular}{c}
     $\boldsymbol{F}^m$ \\
     $-\tau\boldsymbol{G}^m$ 
\end{tabular}
\right], \quad A' =\left[\begin{tabular}{cc}
     $0$& $0$ \\
    $B^T$ & $-\frac{\beta}{K} C$ 
\end{tabular} \right] .
$$

\section{The parallel-in-time method based on preconditioner}\label{sec:pitbP}
\subsection{Design of preconditioners}\label{preconditioner}

Biot's model is a typical saddle-point problem which has eigenvalues equal to $0$ and less than $0$. In order to ensure a stable convergence in solving the discrete linear system with the parallel-in-time method. We design two preconditioners for the discret Biot's model \eqref{eq:fullysystem}: 
\begin{equation}
    \widetilde{P}_1 = 
\left[
\begin{tabular}{c c}
     $A$ & $0$ \\
      $B^{T}$ & $-\widetilde{S}_p$ 
\end{tabular}
\right],
\label{p1}
\end{equation}
and
\begin{equation}
    \widetilde{P}_2 = 
\left[
\begin{tabular}{c c}
     $A$ & $0$ \\
      $0$ & $-\widetilde{S}_p$ 
\end{tabular}
\right],
\label{p2}
\end{equation}
where $\widetilde{S}_p$ is the approximate Schur complement matrix,
\begin{equation}
\widetilde{S}_p =  (\tau  + \frac{\beta}{K})C.
\end{equation}

In Section \ref{chapter_result}, we  numerically investigate the stability of convergence in solving discrete Biot's model \eqref{eq:fullysystem} with precondiitoners $\widetilde{P}_1$ and $\widetilde{P}_2$, as well as their convergence efficiency.

\subsection{Inverted time-stepping method and parallelization}\label{inversed_section}
\begin{figure}
    \centering
    \begin{subfigure}[t]{1.0\textwidth}
           \centering
           \includegraphics[width=\textwidth]{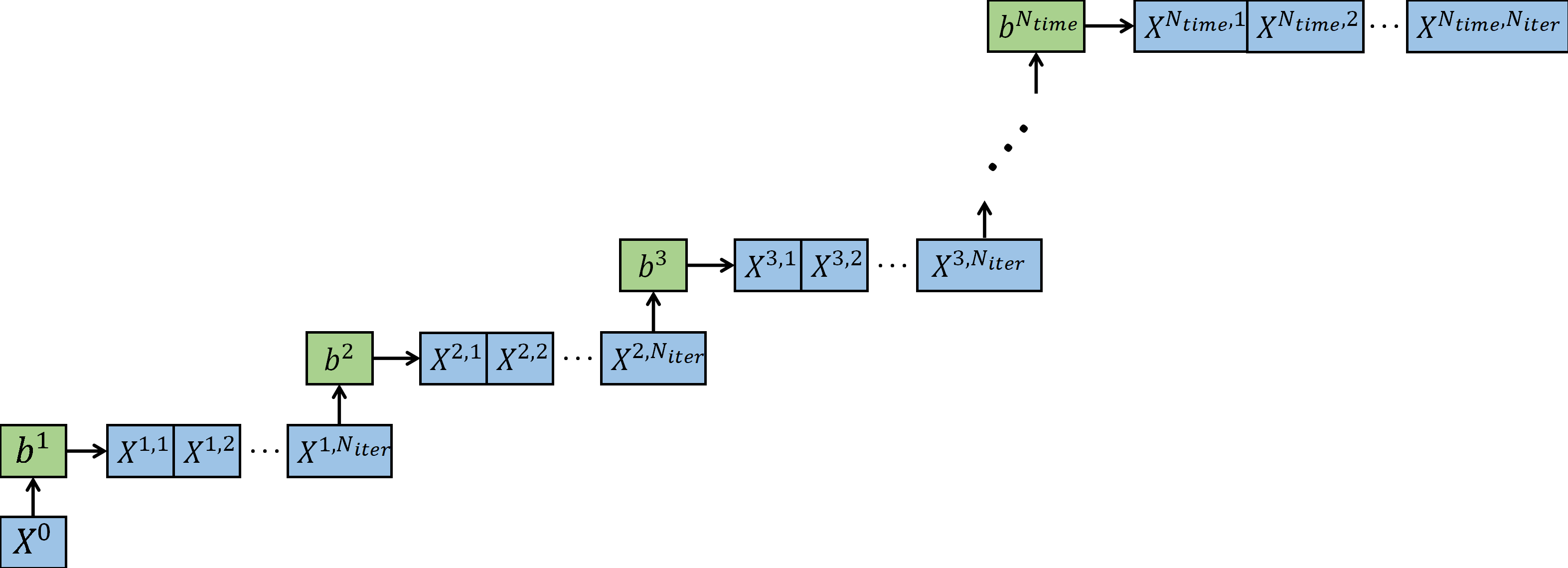}
            \caption{sequential implicit time stepping method}
    \end{subfigure}
    \begin{subfigure}[t]{1.0\textwidth}
           \centering
           \includegraphics[width=\textwidth]{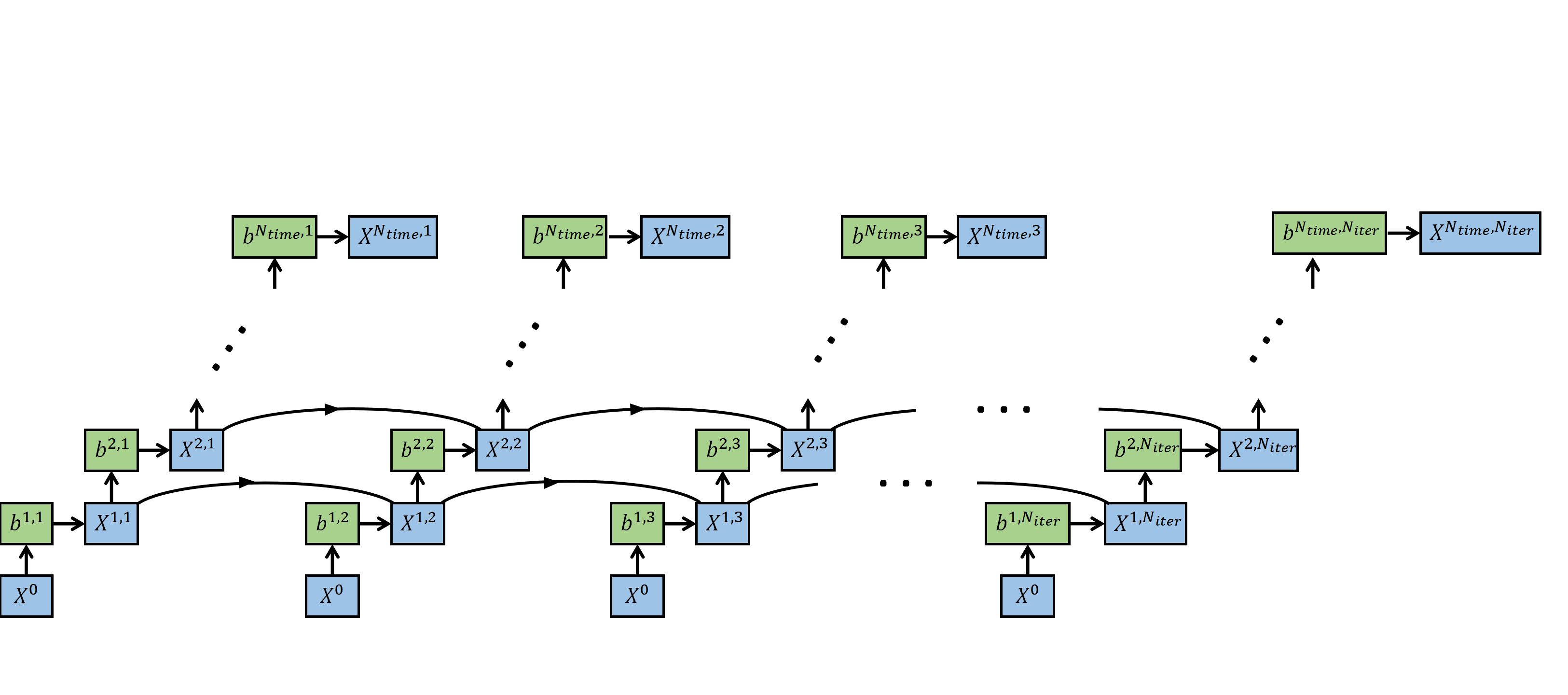}
            \caption{inverted time stepping method}
    \end{subfigure}

    \caption{Comparison between sequential implicit time stepping method and inverted time stepping method.}
    \label{fig:compare}
\end{figure}

The sequential implicit time stepping method, shown in Algorithm \ref{ag:seq}, can not run in parallel along the temporal dimension as it requires the iteration solver for linear system $\mathbbm{A}X^{n}=b$ to converge fully before updating the unknown vector at the next time step $X^{n+1}$. Hence, we propose an inverted time stepping method, shown in Algorithm \ref{inver_normal}. In the inverted time stepping method, the loop of the time steps is taken as the inner loop while the iteration to solve the linear system is taken as the outer loop. A comparison between the sequential implicit time-stepping method and the inverted time stepping method is illustrated in Fig. \ref{fig:compare}.

\begin{algorithm}
\setstretch{1.2}
\begin{algorithmic}[1]
    \caption{sequential implicit time-stepping method}\label{ag:seq}
    \For{$n=1:N_{time}$}
        \State$b^{n} = A'X^{n-1,N_{iter}} + f^{n}$\;
        \State$X^{n,0} =0$ \Comment{zero initial guess for iterative solver}
        \For{$m =1:N_{iter}$}
            \State$X^{n,m}=X^{n,m-1}+P^{-1}(b^{n}-\mathbbm{A}X^{n,m-1})$\;
        \EndFor
    \EndFor
\end{algorithmic} 
\end{algorithm}

\begin{algorithm}
\setstretch{1.2}
\begin{algorithmic}[1]
    \caption{inverted time-stepping method}\label{inver_normal}
    \State $X^{1:N_{time},0} =0$ \Comment{zero initial guess for iterative solver at every time step}
    \For{$m =1:N_{iter}$}
        \For{$n=1:N_{time}$}
         \State $b^{n,m} = A'X^{n-1,m} + f^{n}$\;
         \State  $X^{n,m}=X^{n,m-1}+P^{-1}(b^{n,m}-\mathbbm{A}X^{n,m-1})$\;
        \EndFor
    \EndFor
\end{algorithmic} 
\end{algorithm}

\begin{figure}
    \centering
    \includegraphics[width = 0.75\linewidth]{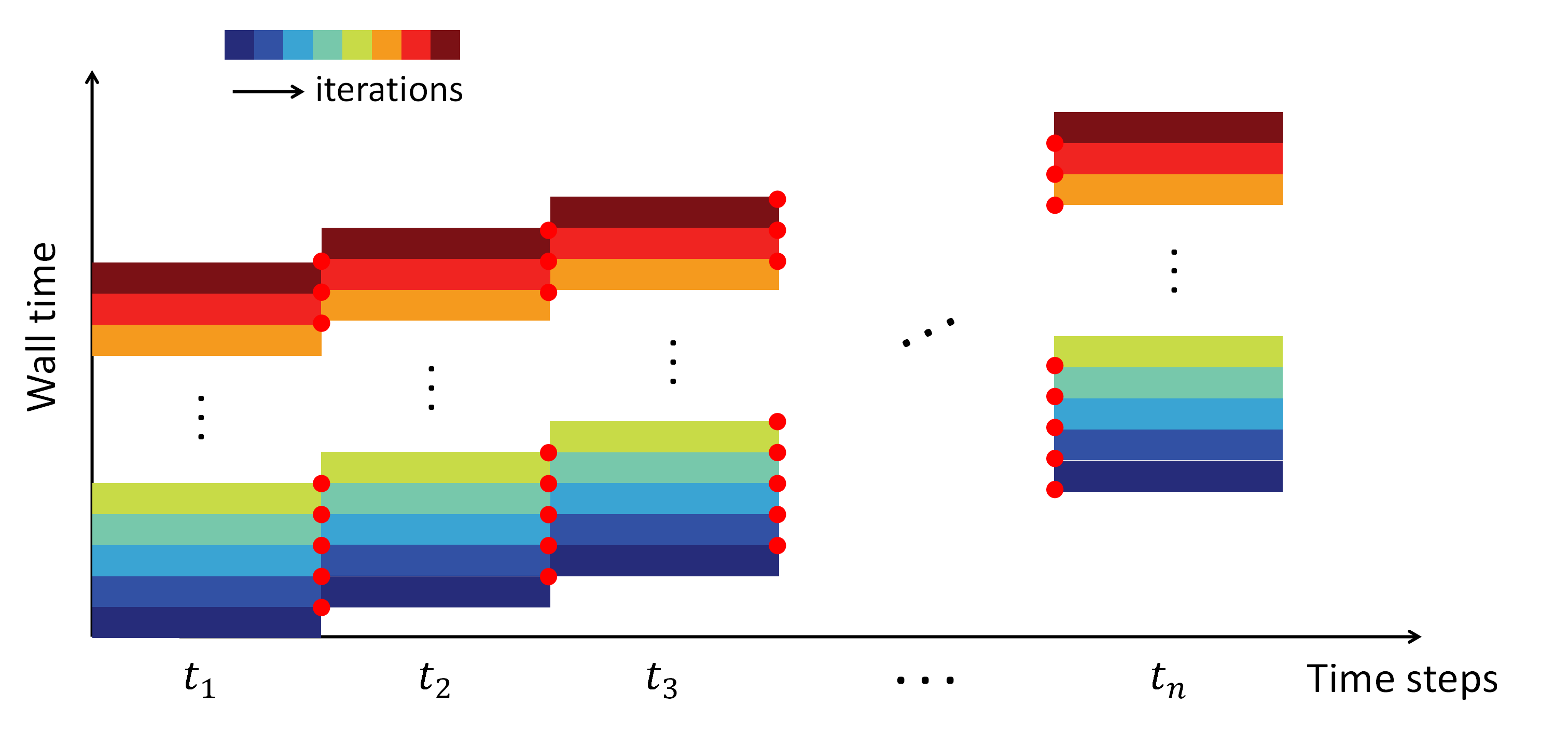}
    \caption{Parallel-in-time Method. The rectangles in the same color stand for iterations calculated by the same thread and red nodes stand for information exchange between two threads.}
    \label{fig:parallel}
\end{figure}

In the inverted time-stepping method, every iteration requires the solution $X^{n, m-1}$ from the last iteration to construct the right-hand vector $b^{n,m}=A'X^{n-1,m}+f^{n}$; it also depends on the solution $X^{n-1,m} $ from last time step. So it is possible to apply the parallelization to the outer loop of linear solver iterations $m$. Fig. \ref{fig:parallel} demonstrates the inverted time-stepping method running in parallel, which we name the parallel-in-time method. The technical details of thread management and data exchange will be discussed in Section \ref{parallelsection}.

One critical difference between the traditional method and the inverted time-stepping method is that the traditional method uses a constant right-hand vector $b^{n}$ for all iterations at one time step $t^n$. In contrast, the inverted method updates the right-hand vector $b^{n,m}$ before every iteration at every time step. So the solution $X^{n,m}$ in the two methods is not exactly equivalent. Proof of the consistency and convergence of the inverted method is under study, which is meaningful to extend the method to more general application scenarios. This paper only demonstrates the feasibility of the inverted method with numerical tests in Section \ref{chapter_result}, and the theoretical proof is beyond the paper's scope.

In practice, the total thread number is not necessary to equal to the number of iterations. Define the function $\mathcal{K}(P,\mathbbm{A},X,b,n_{iter})$ which stands for iterations solver to solve $\mathbbm{A}X=b$ with preconditioner $P$ and initial guess $X$; $n_{iter}$ is the number of iterations executed for each call to function $\mathcal{K}$. The inverted time-stepping method can be written as Algorithm \ref{inver_define}, where the total number of iterations is $N_{iter}$ and the number of threads $N_{thread} =\frac{N_{iter}}{n_{iter}}$ if the inverted method runs in parallel. It can be seen that the inverted method is exactly consistent with the backward Euler time-stepping method when $N_{thread} =\frac{N_{iter}}{n_{iter}}=1$.

\begin{algorithm}
\setstretch{1.2}
\begin{algorithmic}[1]
    \caption{inverted time-stepping method with arbitrary thread numbers}\label{inver_define}
    \State $X^{1:N_{time},0} =0$ \Comment{zero initial guess for iterative solver at every time step}
    \For{$m =1:\frac{N_{iter}}{n_{iter}}$}
        \For{$n=1:N_{time}$}
         \State $b^{n,m} = A'X^{n-1,m} + f^{n}$\;
            \State $X^{n,m}=\mathcal{K}(P,\mathbbm{A},X^{n,m-1},b^{n,m},n_{iter})$\;
        \EndFor
    \EndFor

\end{algorithmic}
\end{algorithm}

Since $X^{n,m+1}$ depends on $X^{n,m}$, thread $m$ can only conduct $X^{n,m+1}=X^{n,m}+P^{-1}(b^{n}-X^{n,m})$ after thread $m-1$ conducts $X^{n,m}=X^{n,m-1}+P^{-1}(b^{n}-X^{n,m-1})$. It means thread $m$ lags behind thread $m-1$ by one time step, and therefore begins and ends one step later than thread $m-1$. Extending this relationship to all threads, it is easy to find that the last thread ends $N_{thread} -1$ time steps later than the first thread. As the workload of every thread is assumed uniform, it leads to a theoretical upper limit for parallel efficiency :
\begin{equation}\label{e1_theo}
    E_1 =\frac{N_{time}}{N_{thread}\cdot\frac{N_{time}+N_{thread}-1}{N_{thread}}}= \frac{N_{time}}{N_{thread}+N_{time}-1}.
\end{equation}
  
This feature leads to the parallel efficiency less than 1. However, as $\frac{N_{time}}{N_{thread}}$ increases, the theoretical parallel efficiency will infinitely close to 1, which means the parallel strategy is more efficient for cases with a higher number of time steps.

For the traditional implicit time-stepping method, right-hand vector $b^{n}$ only needs to be computed once at time step $n$, however, for Algorithm \ref{inver_define}, $b^{n,m}$ is calculated by each thread $m$ at each time step $n$. It means the parallel-in-time method only accelerates the iteration process while the wall-clock time costed by computing $b$ is constant. Considering the computing time for $b^{n,m}$, it leads to another theoretical upper limit for efficiency:
\begin{equation}\label{E2_theo}
         E_2 = \frac{T_b +N_{iter}\cdot t_{iter}}{T_b +n_{iter}\cdot t_{iter}} \cdot \frac{1}{N_{thread}} = \frac{\frac{T_b}{N_{iter}}+t_{iter}}{\frac{N_{thread}}{N_{iter}}\cdot T_b + t_{iter}},
\end{equation}
where $t_{iter}$ is the computing time for one Krylov Subspace iteration and $T_b$ is the time for computing $b$.
     
This feature also leads to the parallel efficiency less than 1, while as $\frac{N_{iter}}{N_{thread}}$ increases, the theoretical parallel efficiency will infinitely close to 1. Generally speaking, the larger the degrees of freedom of a discrete system, the greater $N_{iter}$ required to solve the linear system of equations.

Considering two theroretical upper limits \eqref{e1_theo} and \eqref{E2_theo}, the total theoretical parallel efficiency can be derived as $E$:
\begin{equation}
         E =E_1 \cdot E_2.
\end{equation}
     
For certain number of threads, setting more time steps results in larger $\frac{N_{time}}{N_{thread}}$  thus $E_1$ is closer to 1; using more spatial steps results in larger $\frac{N_{iter}}{N_{thread}}$ so $E_2$ is closer to 1. This means that the parallel-in-time method proposed in this paper has better theoretical parallel efficiency for larger-scale problems.

\subsection{Implementation of the parallel-in-time method}\label{parallelsection}
The inverted time-stepping method has been introduced in Section \ref{inversed_section}, briefly explaining the parallel strategy. This section will illustrate the details of threads management and data exchange.

Firstly, the pure shared memory parallel architecture implemented by OpenMP is illustrated. Although threads in the inverted method can run in parallel, they do not start and end simultaneously. So the outer loop of $m$ cannot be parallelized simply by clause “\#pragma omp for” which starts the thread for every iteration simultaneously; instead, the parallel method should start threads in turn as Algorithm \ref{OpMP}.

\begin{algorithm}
\setstretch{1.2}
\begin{algorithmic}[1]
    \caption{parallel-in-time method with OpenMP}
    \label{OpMP}
   \State $X^{1:N_{time},0} =0$ \Comment{zero initial guess for iterative solver at every time step}
    \State $\text{set\_thread\_num}(N_{thread})$\;
     \State \#pragma omp parallel\;
        \For{$n=1:N_{time}+N_{thread} -1$}
            \State $n_t = n- TID$\;
        \If{$n_t>0\ \   \text{and}\ \  n_t<N_{time}+1$}
          \State $b^{n_t} = A'X^{n_t-1} + f^{n_t}$\;
            
            \State $X^{n_t}=\mathcal{K}(P,\mathbbm{A},X^{n_t},b^{n_t},n_{iter})$\;
            \EndIf 
        \State \#omp barrier
        \EndFor
\end{algorithmic}
\end{algorithm}

With the shared memory architecture, solutions $X$ is a 2-dimensional array of size  $N_{time}\times N_{DOF}$. At each time step, $X^{n_t}$ stands for the solution updated by the latest iteration, so it is unnecessary to extend another dimension to store the solutions output by earlier iterations. In practice, array length along the dimension of time can be shorten to $N_{thread}$ from $N_{time}$, because the solution $X^{n}$ will no longer be used in iterations after being updated by the last thread $N_{thread}-1$. Releasing memory in a timely manner can effectively improve computational performance.

In Algorithm \ref{OpMP}, TID is the thread ID starting with 0. Calculating $n_t = n-TID$ and using IF-clause $n_t>0\quad  \text{and}\quad n_t<N_{time}+1$ implements the asynchronous starting and ending of threads illustrated in Fig. \ref{fig:parallel}.  For each thread, $n_t$ is the time step which is going to be computed. For example, at time step $n = 1$, thread $0$ has $n_t =1$ and is going to compute $X^1$, and the other thread will not run. Then when $n =2$,  thread $0$ has $n_t=2$, and thread $1$ has $n_t=1$ and is going to update $X^1$. 

Parallel-in-time method is easy to implement with OpenMP in shared-memory multiprocess architecture, which is generally more efficient than point-to-point communication. However, it can only run in one computing nodes which means the maximum threads number is limited by hardware. In order to employ the parallel-in-time method on high performance platforms, a hybrid MPI-OpenMP parallel programs is implemented as Algorithm \ref{hybridmpiOpMP}.

\begin{algorithm}

\setstretch{1.2}
\begin{algorithmic}[1]

    \caption{parallel-in-time method with hybrid MPI-OpenMP programming}
    \label{hybridmpiOpMP}
    \State $X^{1:N_{time}} =0$ \Comment{zero initial guess for iterative solver at every time step}
    \State $\text{set\_process\_num}(N_{process})$\;
    \State $\text{set\_thread\_num}(N_{thread})$\;
    \State \#pragma omp parallel\;
        \For{$n=1:N_{time}+N_{thread}\cdot N_{process} -1$}
        \State $n_t = n - PID\cdot N_{thread}- TID$\;
        \If{$TID =0 \ \ \text{and}\ \  mod(PID,2)=0$}
            \If{ $n_t -N_{thread} >0 \ \ \text{and} \ \ n_t -N_{thread}<N_{time}+1 \ \ \text{and} \ \ PID != N_{process}-1$ }
            \State MPI\_Send($X^{n_t-N_{thread}}$, destination=$PID+1$, tag=$n_t-N_{thread}$)\;
            \EndIf
            \If{$n_t >0 \ \ \text{and} \ \ n_t <N_{time}+1 \ \ \text{and} \ \ PID != 0$ }
            \State MPI\_Recv($X^{n_t}$, destination=$PID-1$, tag=$n_t$)\;
            \EndIf
        \EndIf
        \If{$TID =0 \ \ \text{and}\ \  mod(PID,2)=1$}
            \If{$n_t >0 \ \ \text{and} \ \ n_t <N_{time}+1 \ \ \text{and} \ \ PID != 0$ }
            \State MPI\_Recv($X^{n_t}$, destination=$PID-1$, tag=$n_t$)\;
            \EndIf
            \If{ $n_t -N_{thread} >0 \ \ \text{and} \ \ n_t -N_{thread}<N_{time}+1 \ \ \text{and} \ \ PID != N_{process}-1$ }
            \State MPI\_Send($X^{n_t-N_{thread}}$, destination=$PID+1$, tag=$n_t-N_{thread}$)\;
            \EndIf
        \EndIf

        \If{$n_t>0\ \ \text{and}\ \ n_t<N_{time}+1$}
         \State $b^{n_t} = A'X^{n_t-1} + f^{n_t}$\;
           \State  $X^{n_t}=\mathcal{K}(P,\mathbbm{A},X^{n_t},b^{n_t},N_{restart},n_{iter})$\;
            \EndIf
        \State \#omp barrier
        \EndFor
        
\end{algorithmic}
\end{algorithm}

In hybrid MPI-OpenMP programming, one process has several threads which shares the memory, and the data exchange between processes relies on MPI. In pseudo code above, $N_{process}$ is the number of processes and $N_{thread}$ is the number of threads in each process, so the total number of threads is $N_{thread}\times N_{process}$ (the numbers of threads are uniform to all process). And $PID$ is the process ID starting from 0 and $TID$ is the thread ID belongs to certain process.

\clearpage
\section{Numerical experiments}
\label{chapter_result}
This section will present several numerical tests of Biot‘s problem computed with the preconditioned parallel-in-time method. These numerical results provide strong evidence that the algorithm achieves the desired performance. The following aspect of performance will be focused on. Firstly, consistency to the exact solution will be illustrated. As discussed in Section \ref{inversed_section}, the numerical system in the parallel-in-time method is not the same as the traditional sequential method, so the numerical results should be proven consistent with the analytical solutions. Then the order of accuracy and the numerical stability of the primal discrete scheme is expected to be maintained. The error of primal spatial discrete scheme runs in sequential can converge by the first order; simultaneously, the numerical oscillation near rapid pressure slope can be eliminated. These two properties are critical for an algorithm to output physical and reliable results and expected in the parallel-in-time method. Finally, the parallel efficiency will be illustrated, as the starting point of designing parallel algorithms is to reduce the wall-clock computing time. For all numerical cases in this section, the solver used is coded in Fortran and the iterative method is GMRES implemented in PETSc.

\subsection{Trigonometric Function Test}
\label{sec:trigono}
Trigonometric functions are commonly used in numerical tests. For these cases, analytical solutions are explicit, and the spatial frequency is simple to analyze. 
Assume the domain is $\Omega=\left[0,1\right]\times\left[0,1\right]$, boundaries are noted as $\Gamma_1 = \{(1,y);0\leq y\leq 1\}$, $\Gamma_2 = \{(x,0);0\leq x \leq 1\}$, $\Gamma_3 = \{(0,y);0\leq y\leq 1\}$ and  $\Gamma_4 = \{(x,1);0\leq x \leq 1\}$. The total time length $T = 10$. The external force term is
\begin{equation}
\begin{aligned}
    \boldsymbol{f} = -3\pi cos(t)\left( 
    \begin{aligned}
          \alpha sin(3\pi x)cos(3\pi y) -3\pi (\lambda + 3\mu)cos(3\pi x)cos( 3\pi y) \\
        \alpha cos(3\pi x)sin(3\pi y) -3\pi (\lambda + 3\mu)cos(3\pi x)cos(3\pi y) 
    \end{aligned}
    \right. \\
    \left. \begin{aligned}
        +3\pi(\lambda +\mu)sin(3\pi x)sin(3\pi y)\\
        +3\pi(\lambda +\mu)sin(3\pi x)sin(3\pi y)
    \end{aligned} \right),
\end{aligned}
\end{equation}
and the source term is
\begin{equation}
g =  3\alpha \pi sin(t)(cos(3\pi x)sin(3\pi y) - sin(3\pi x)cos(3\pi y))+18K\pi^2 cos(t)cos(3\pi x)cos(3\pi y).
\end{equation}
Dirichlet boundary conditions are used at all boundaries

\begin{equation} 
    \left\{
    \begin{aligned}
         p = cos(t)cos(3\pi x)cos(3\pi y)  \\
               u_1 = cos(t)cos(3\pi x)cos(3\pi y) \\
               u_2 = cos(t)cos(3\pi x)cos(3\pi y) 
       \end{aligned}    
    \right.
    \quad  \text{on} \quad \Gamma_i, \ i=1,2,3,4.
\end{equation}
The initial condition is 
\begin{equation}
    \left\{
    \begin{aligned}
        p = cos(3\pi x)cos (3\pi y)\\
        u_1 =cos(3\pi x)cos(3\pi y)\\
        u_2 = cos(3\pi x)cos(3\pi y)
    \end{aligned}
    \right. .
\end{equation}
The exact solution can be easily derived from above conditions
\begin{equation}
    \left\{
    \begin{aligned}
        p = cos(t)cos(3\pi x)cos(3\pi y)  \\
               u_1 = cos(t)cos(3\pi x)cos(3\pi y) \\
               u_2 = cos(t)cos(3\pi x)cos(3\pi y) 
    \end{aligned}
    \right.
    \quad \text{in}\quad  \Omega .
\end{equation}
The parameters are set as following:
$$
E = 1, \nu =0.3, \alpha =1, K = 1,
$$
Where $E$ is the Young's modulus, $\nu$ is the Poisson's ratio, $\alpha$ is the Biot-Wills constant, and $K$ is the ratio of material
permeability to fluid viscosity. 
The domain is discretized with uniform Triangular P1-P1 elements, and the spatial step length of the mesh is $h = 1/2^k$, where k is an integer. In the mesh refinement process, the new mesh is obtained by connecting the midpoints of each side of the triangle, i.e., the original triangular element is divided into four triangular elements after one refinement. The Krylov subspace method used in this case is GMRES 
with preconditioner $\widetilde{P_1}$ in Eq. \eqref{p1}. To set the numbers of iterations properly for different methods, the tolerance of the final residual in iterations is set the same as $10^{-7}$ for all cases, which is small enough to consider the linear equations system $\mathbbm{A}X=b$ solved accurately.

\begin{figure}
    \begin{subfigure}[t]{0.5\textwidth}
           \centering
           \includegraphics[width=\textwidth]{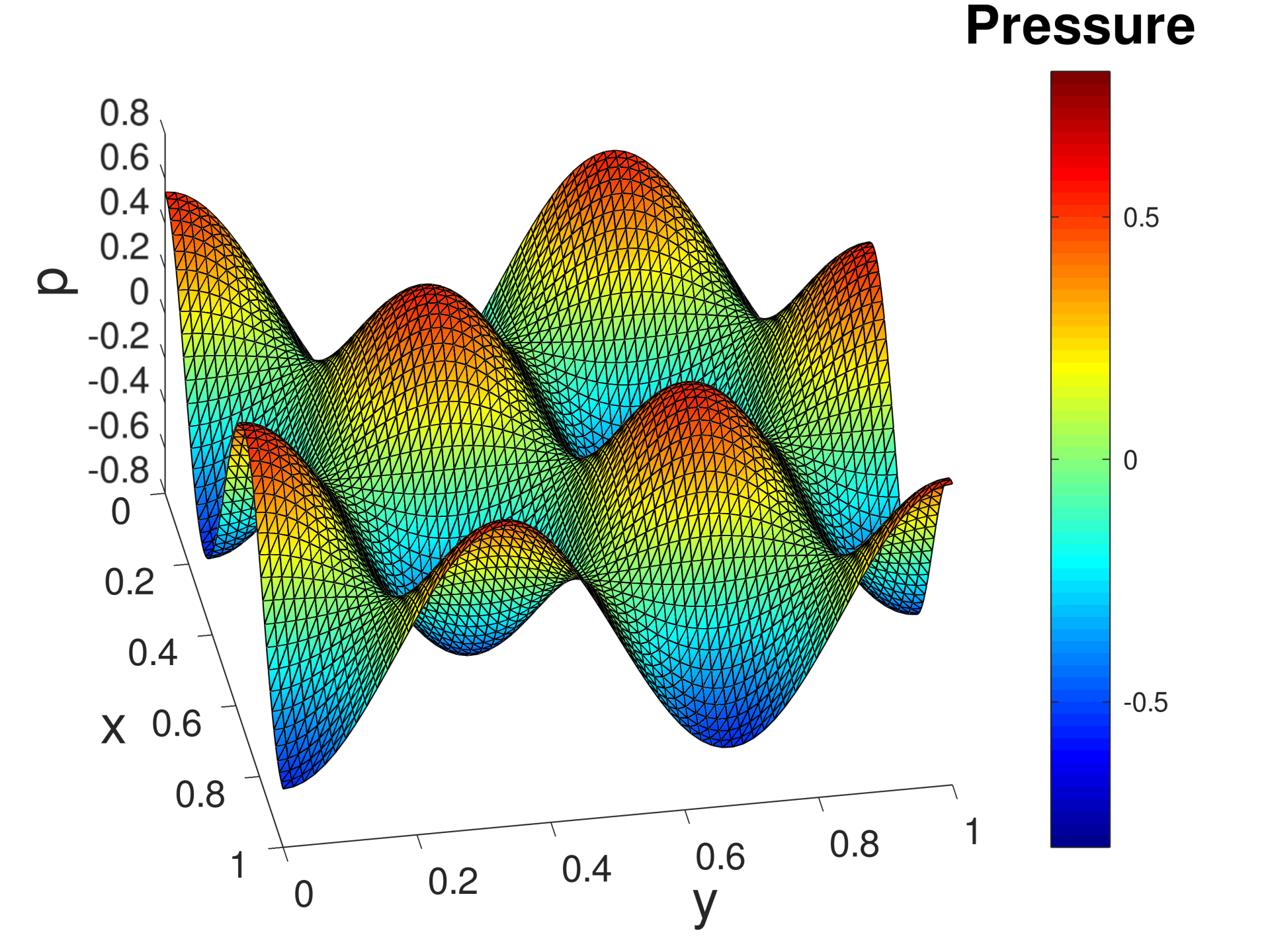}
            \caption{sequential (solution)}
            \label{fig:Case1}
    \end{subfigure}
    \begin{subfigure}[t]{0.5\textwidth}
            \centering
            \includegraphics[width=\textwidth]{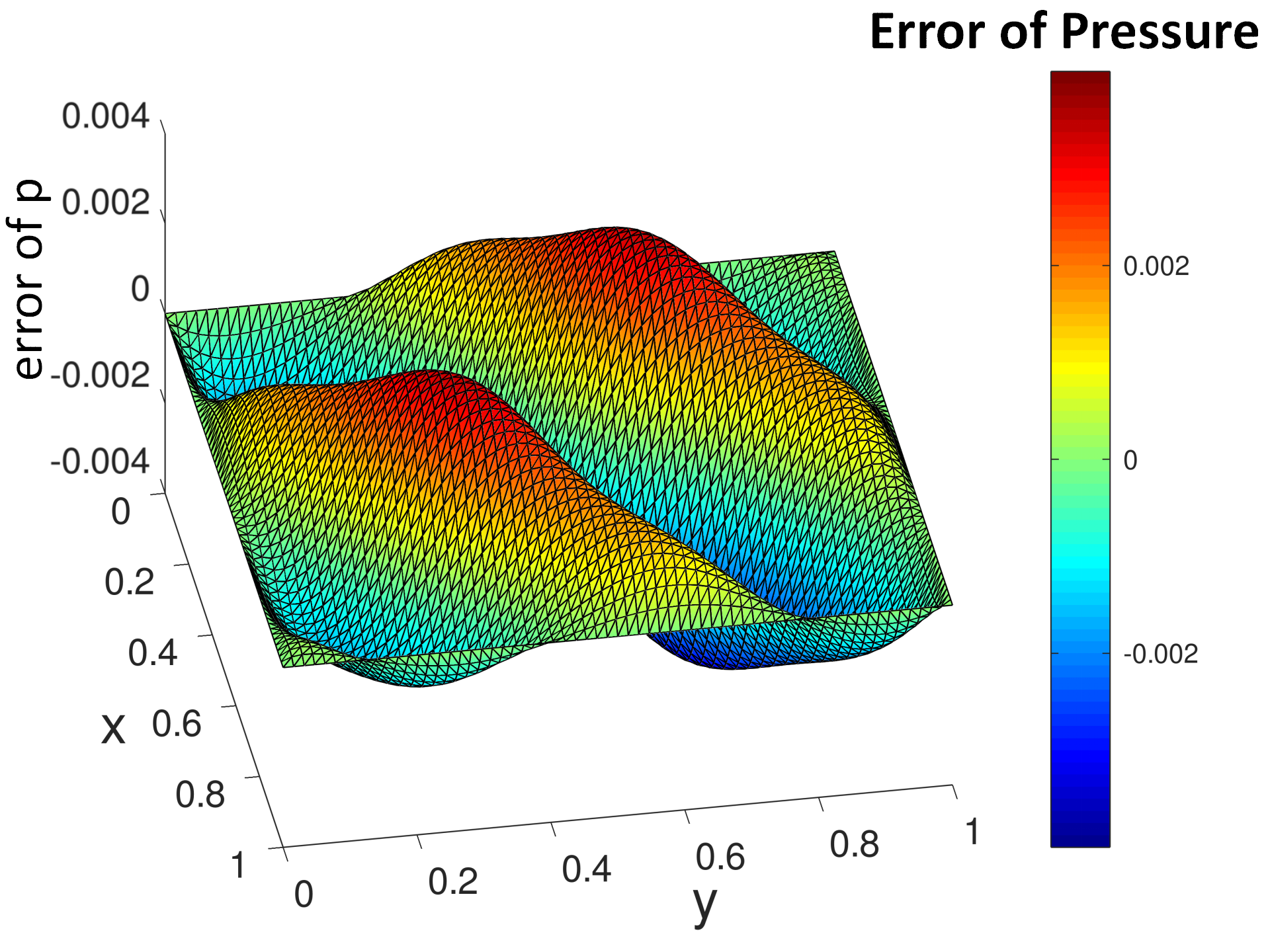}
            \caption{sequential (error)}
            \label{fig:b}
    \end{subfigure}
\begin{subfigure}[t]{0.5\textwidth}
           \centering
           \includegraphics[width=\textwidth]{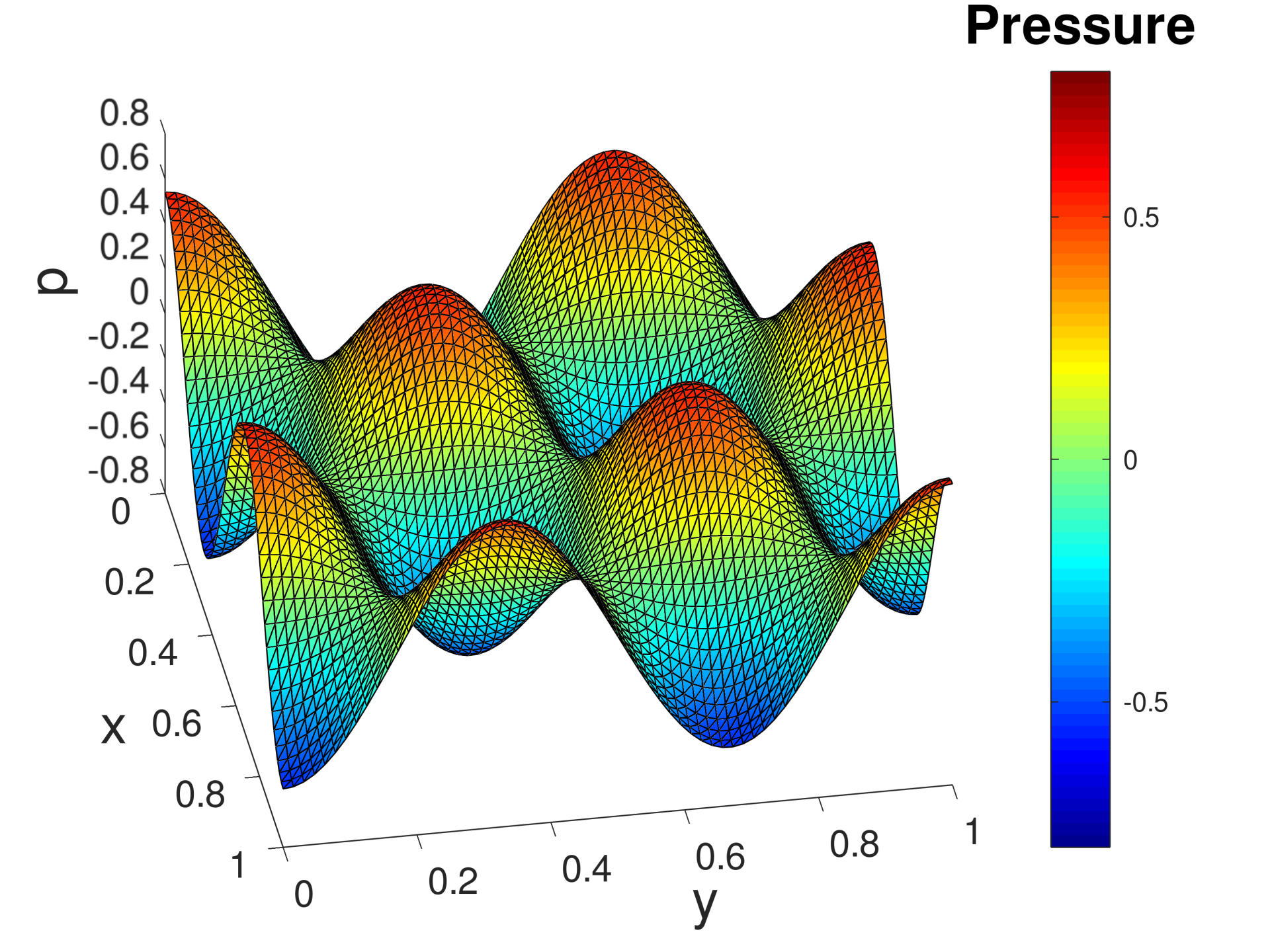}
            \caption{parallel 16 threads (solution)}
            \label{fig:Case1}
    \end{subfigure}
    \begin{subfigure}[t]{0.5\textwidth}
            \centering
            \includegraphics[width=\textwidth]{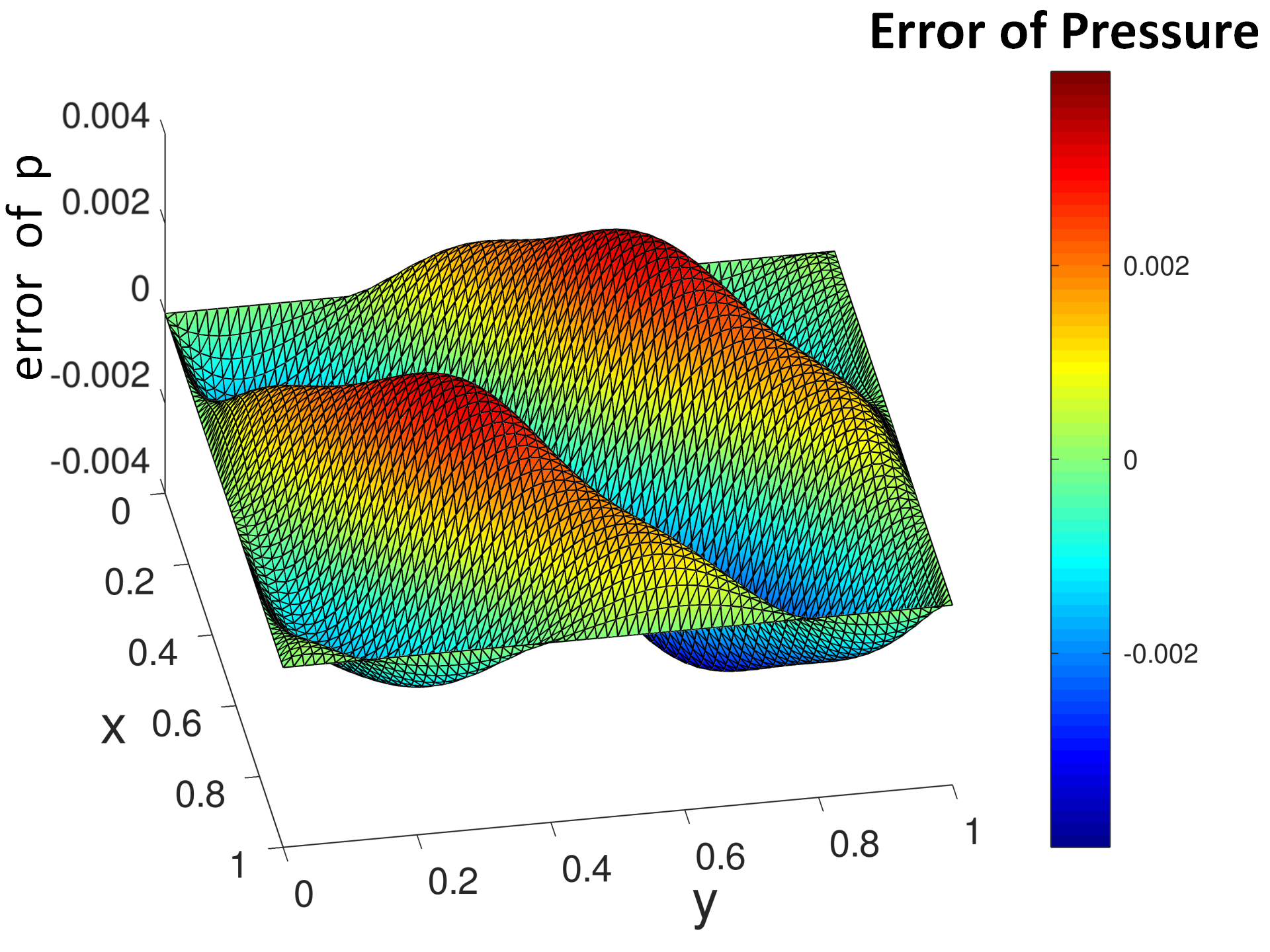}
            \caption{parallel 16 threads (error)}
            \label{fig:b}
    \end{subfigure}
    \begin{subfigure}[t]{0.5\textwidth}
           \centering
           \includegraphics[width=\textwidth]{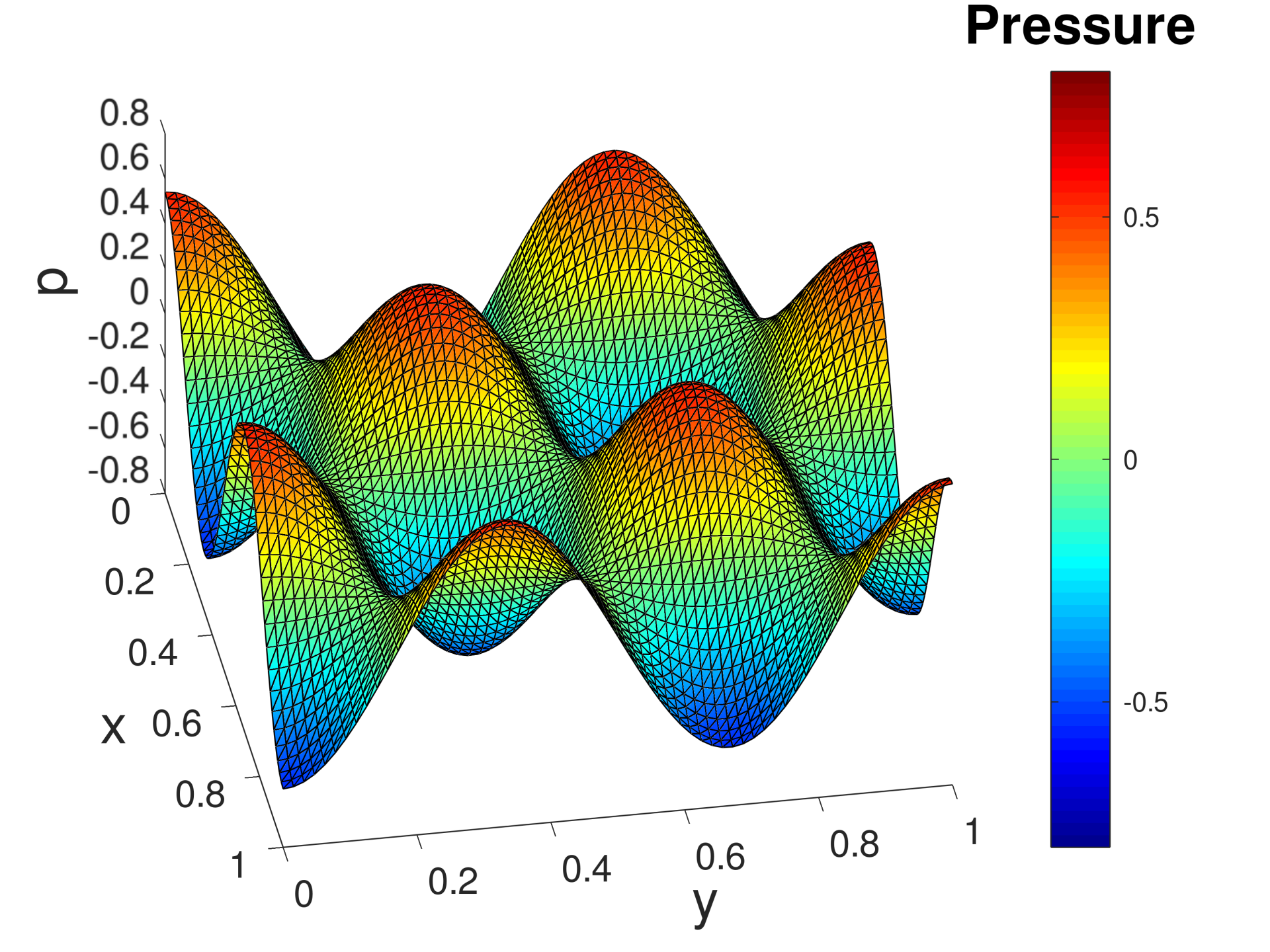}
            \caption{parallel 64 threads (solution)}
            \label{fig:Case1}
    \end{subfigure}
    \begin{subfigure}[t]{0.5\textwidth}
            \centering
            \includegraphics[width=\textwidth]{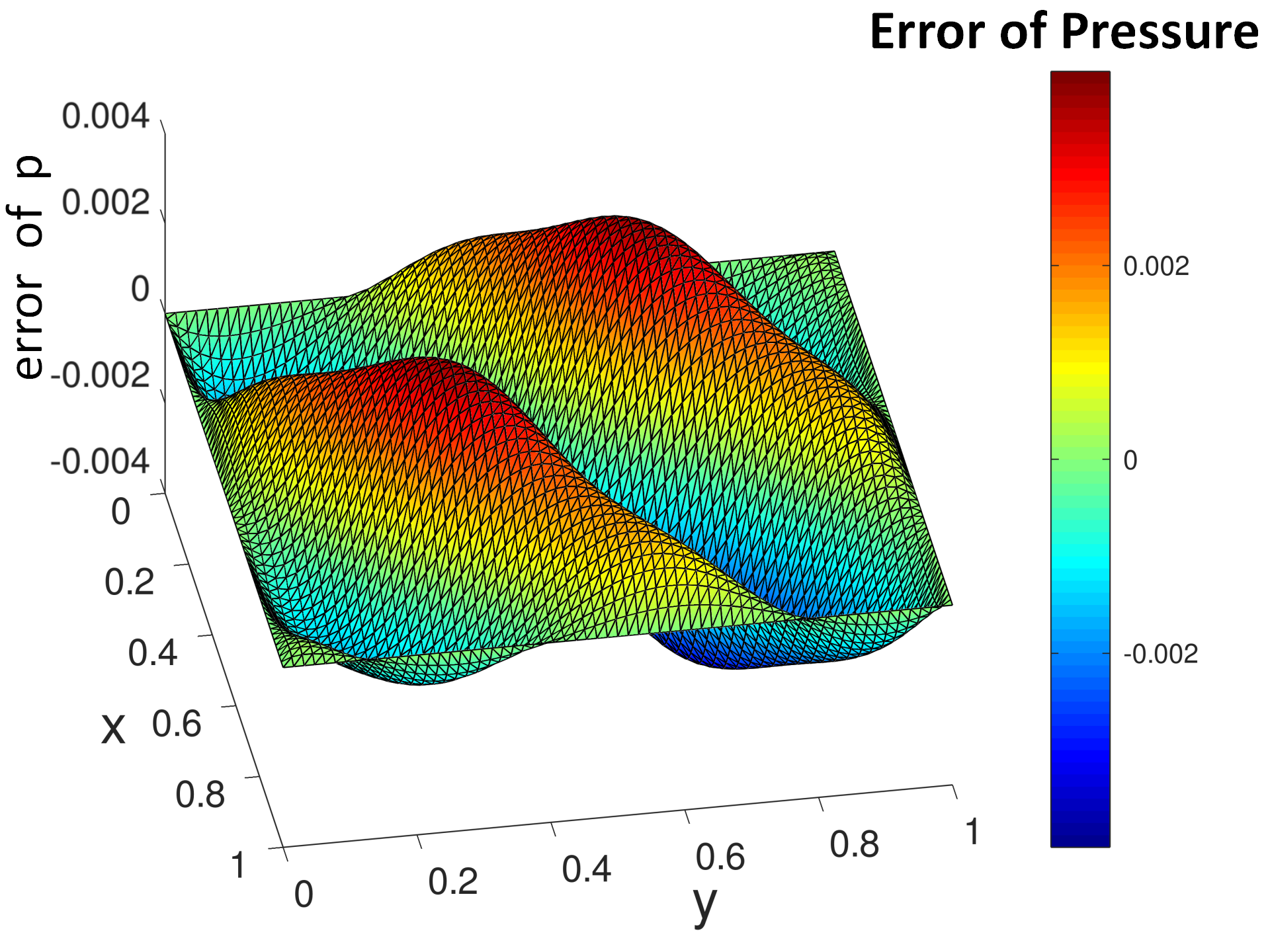}
            \caption{parallel 64 threads (error)}
            \label{fig:b}
    \end{subfigure}

    \caption{Comparison of solution and error of $p$ at $t=1$ with backward Euler method and the patallel-in-time method.}
    \label{fig:pandperror}
\end{figure}

Choose spatial step length $h =1/64$ time step length $\tau =2^{-10}$ and the solutions can be found in 
Fig. \ref{fig:pandperror}. The results numerically prove the consistency of the parallel-in-time method. Though, the numerical system of the parallel-in-time method is different from the sequential time-stepping method, the parallel-in-time method will not introduce additional systematical error compared with the sequential method as long as the final residual of iterations is small enough.

\begin{figure}

    \begin{subfigure}[t]{0.5\textwidth}
           \centering
           \includegraphics[width=1.0\textwidth]{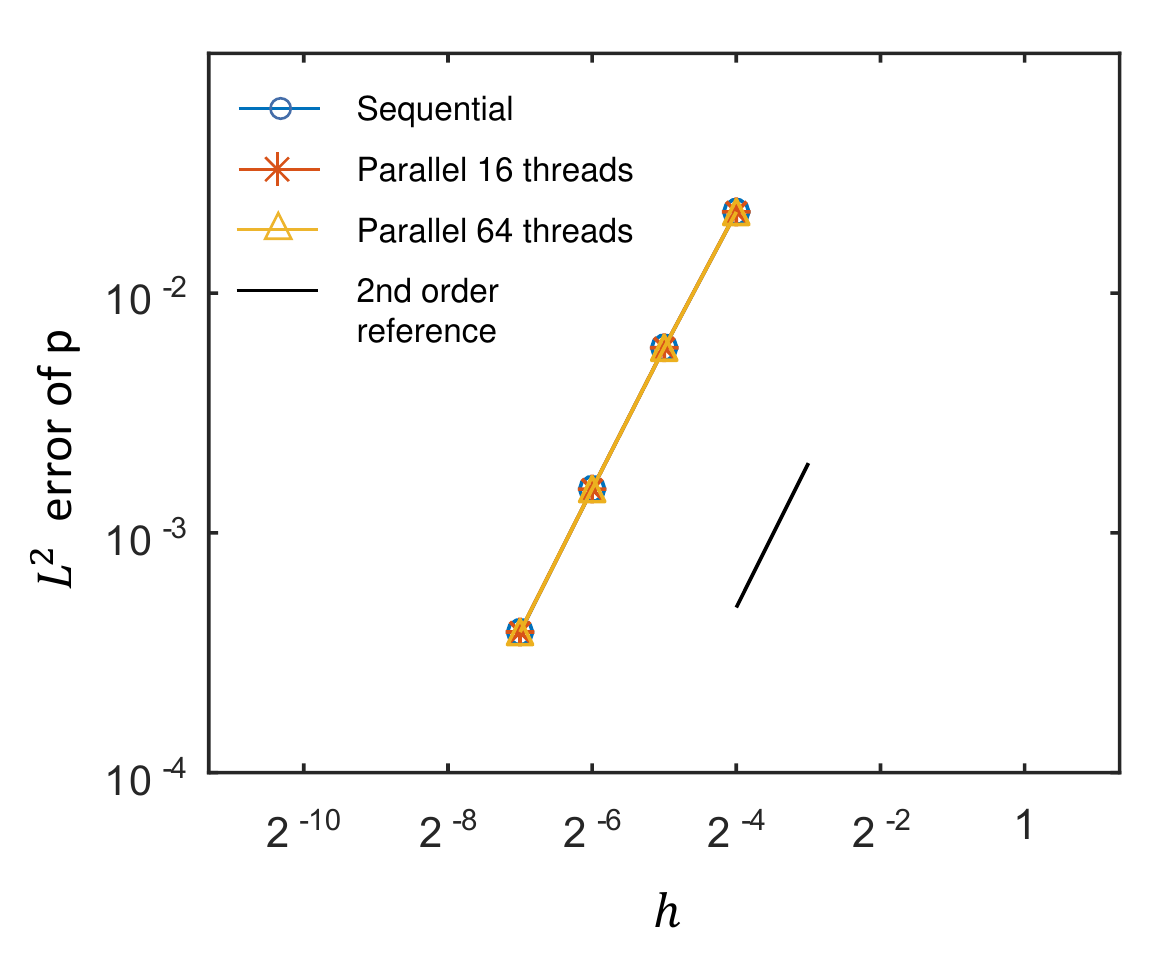}
           \caption{spatial order of accuracy}
    \end{subfigure}
    \begin{subfigure}[t]{0.5\textwidth}
            \centering
            \includegraphics[width=1.0\textwidth]{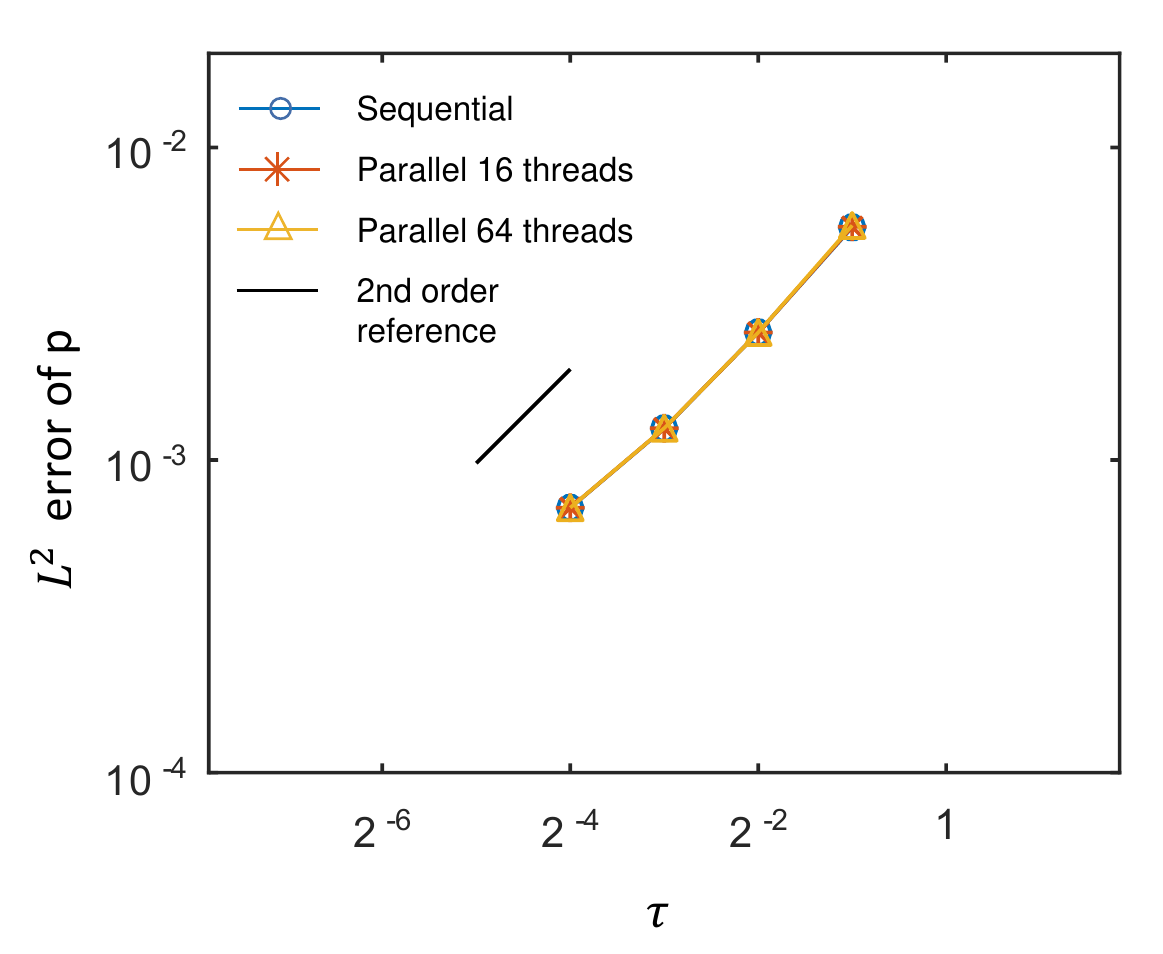}
            \caption{temporal ordder of accuracy}
    \end{subfigure}

    \caption{Convergence rate of $L^2$ error of pressure field.}
    \label{fig:cos_converge}
\end{figure}

To test the convergence rates of error with spatial and temporal discrete step lengths, choose $\tau =2^{-10}$,  $h=h_0/2^k$ where $h_0 =1/16$ and $k =0,1,2,3$; and choose $h =1/128$, $\tau =\tau_0/2^k$ where $\tau_0 =1/2$ and $k = 0,1,2,3$. The results are presented in Table \ref{cos_errortable} and Fig. \ref{fig:cos_converge}. The results show that the parallel-in-time method maintains the order of accuracy of the primal discrete scheme in both time and space dimensions. The parallel-in-time method's error converges in the same order as the error computed by the sequential backward Euler time-stepping method. In general, stabilized P1-P1 scheme has first-order accuracy, while in this case, a second-order convergence is observed when decreasing the grid size from $1/16$ to $1/128$. This superconvergence phenomenon is due to the fact that the problem has a smooth solution and the grids are uniform in these cases; see the details in \cite{superconverg}. 

\begin{table}
    \centering
    \caption{$L^2$ errors of $p$.}
    \label{cos_errortable}
    \begin{tabular}[width = \textwidth]{c c c c}
    \multicolumn{2}{c}{$\tau =1/1024, \quad t =1$}\\ \hline
      method  & Sequential & Parallel 16 threads & Parallel 64 threads \\ \hline
      $h$ & $L^2$ error\\
      1/16   & 0.02180138 &0.02180176 &0.02180133 \\ 
      1/32  & 0.00592069 & 0.00592098 &0.00592173 \\
      1/64 & 0.00152259  & 0.00152273 &0.00152251\\
      1/128 & 0.00038509 & 0.00038521 &0.00038576\\ \hline
    \end{tabular}
    \begin{tabular}[width = \textwidth]{c}
        \quad   \\
    \end{tabular}
    
      \begin{tabular}[width = \textwidth]{c c c c}
    \multicolumn{2}{c}{$h=1/128, \quad t =1$}\\ \hline
      method  & Sequential & Parallel 16 threads & Parallel 64 threads \\ \hline
      $\tau$ & $L^2$ error \\
      1/2   & 0.00556048 & 0.00556423 &0.00556282\\ 
      1/4  & 0.00255742 & 0.00255810 &0.00255919 \\
      1/8 & 0.00126302  & 0.00126384 &0.00126470\\
      1/16 & 0.00070293 & 0.00070325 &0.00070382\\ \hline
    \end{tabular}
    
\end{table}

\begin{figure}

    \begin{subfigure}[t]{0.5\textwidth}
           \centering
           \includegraphics[width=1.1\textwidth]{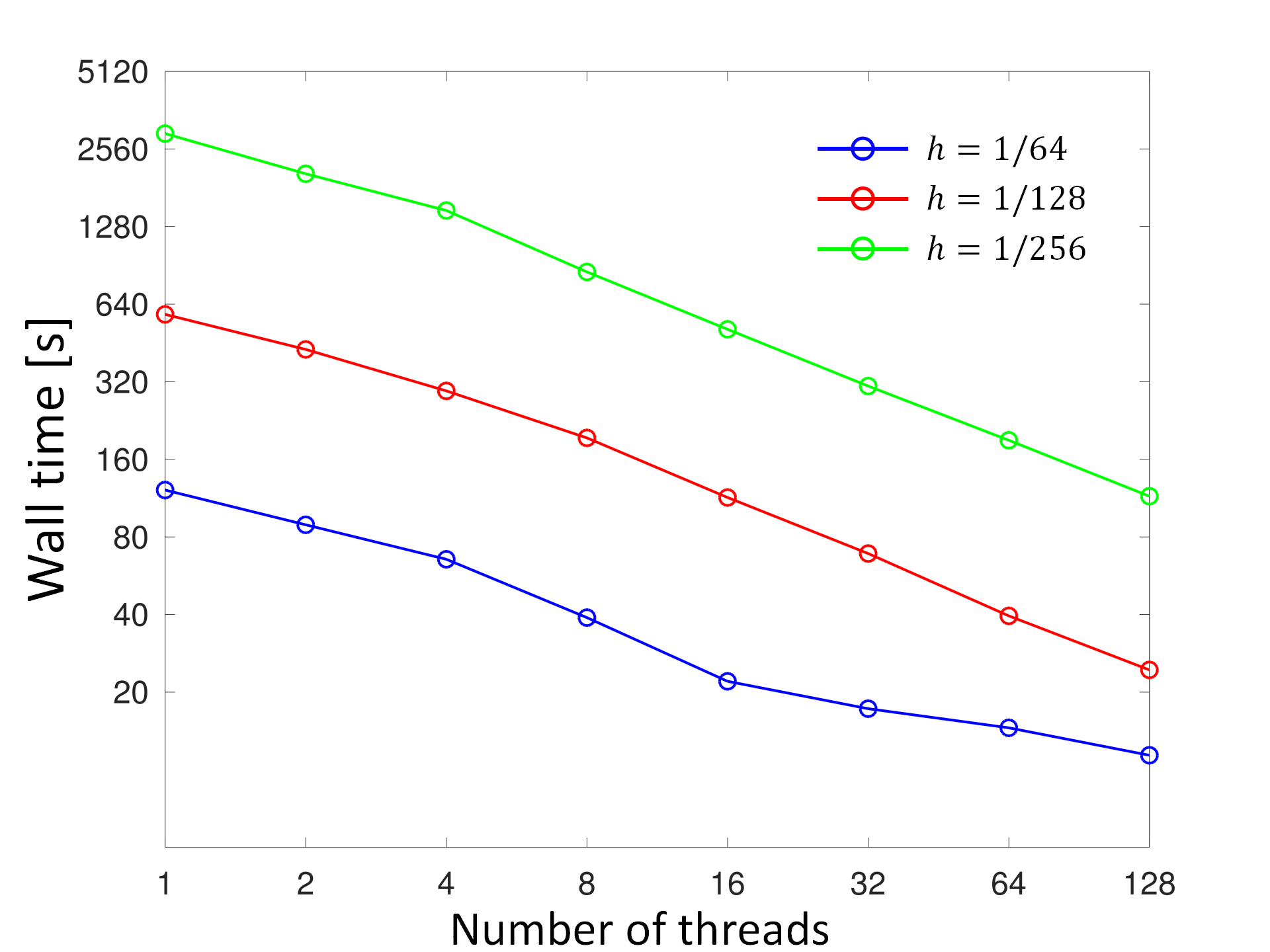}
            \caption{Wall time with different spatial step sizes,\\ $\tau = 1/1024$ }
            \label{fig:Case1}
    \end{subfigure}
    \begin{subfigure}[t]{0.5\textwidth}
            \centering
            \includegraphics[width=1.1\textwidth]{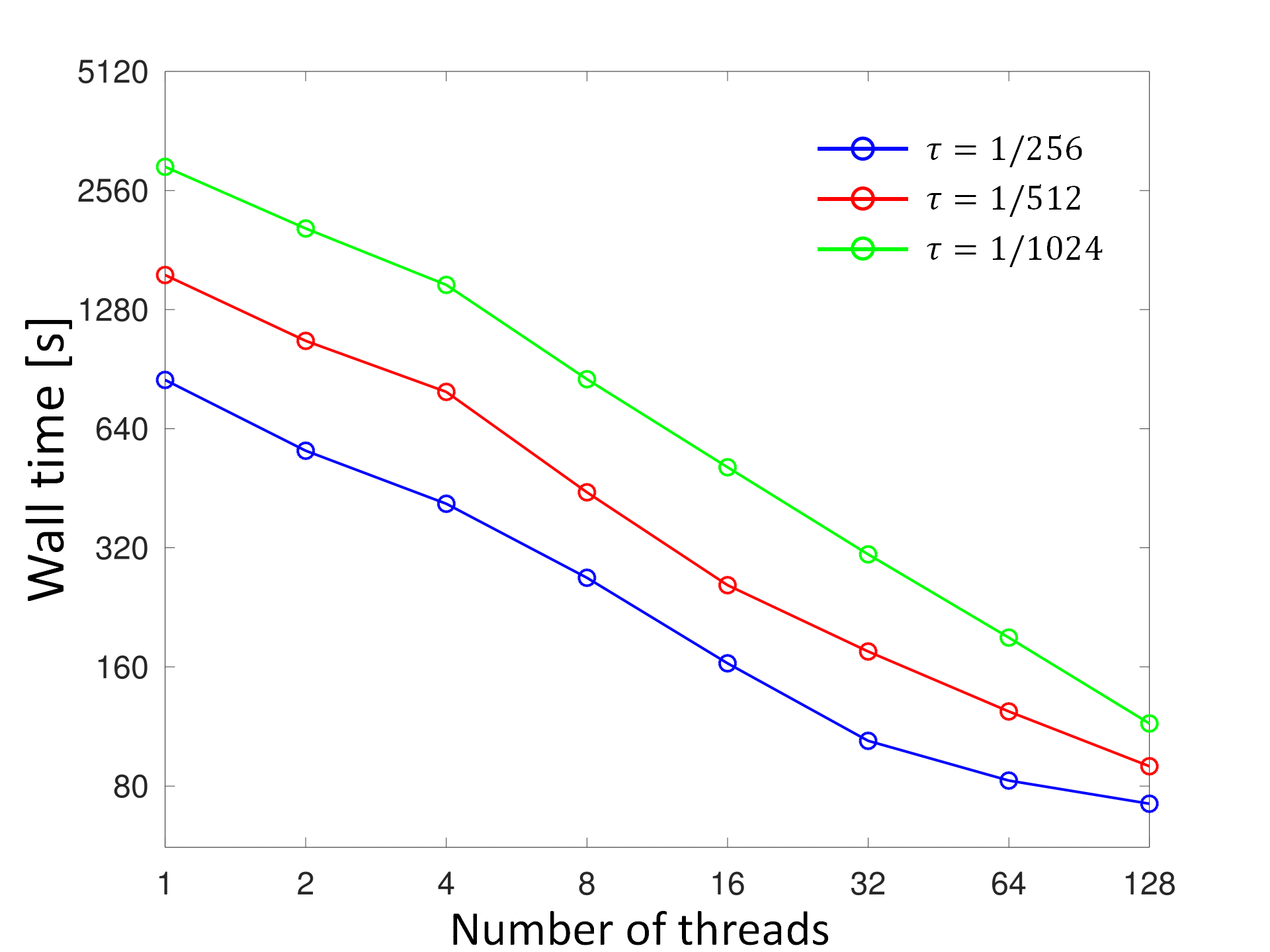}
            \caption{Wall time with different time step sizes,\\ $h =1/256$}
            \label{fig:b}
    \end{subfigure}

    \caption{Wall time of the parallel-in-time method running in the different numbers of threads.}
    \label{fig:coswalltime}
\end{figure}

To study the parallel efficiency, run the parallel-in-time method in $2^k$ threads where $k = 0,1,2...., 7$. As presented in Algorithm \ref{inver_define}, when the number of threads is chosen as 1, the parallel-in-time method is the same as the sequential backward Euler method. The wall-clock time and speedup ration is shown in Table \ref{coswalltimetable} and Fig. \ref{fig:coswalltime}. Although the parallel efficiency in practice is less than 1, the speedup is still approximately linear at all time as the number of threads increases from 1 to 128. 

Note that when $\tau =1/256$, the parallel efficiency significantly decreases when doubling the number of threads from 64 to 128, as discussed in Section \ref{inversed_section}, the parallel efficiency is negatively related with the ratio of the number of threads to the number of time steps; 
besides, when $\tau =1/1024$ and $h=1/64$, the speedup ratio is worse than cases where $h=1/128$ and $h=1/256$. It can be explained by the property discussed in Section \ref{inversed_section}, for cases with small degree of freedom, the theoretical speedup ratio is significantly less than 1.

Besides, another possible reason that leads to a lower parallel efficiency than the theoretical value is that, with the parallel-in-time method, each thread needs to store a complete solution of all degrees of freedom, which makes the parallel-in-time method requires $N_{thread}$ times as much memory as the sequential method, and therefore the time to read and write data increases significantly.

\begin{table}
    \centering
    \caption{Wall time and speedup ratio of the parallel-in-time method running in the different numbers of threads.}
    \label{coswalltimetable}
    \begin{tabular}[width = \textwidth]{c c c c c c c}
    \multicolumn{1}{c}{$\tau =1/1024, \quad t =1$}\\ \hline
      $h$  & $1/64$ &  $1/128$ & $1/256$ & $1/64$ &  $1/128$ & $1/256$ \\ \hline
      Number of threads  &  Wall time [s]  & & & Speedup ratio       \\
      1   & 121.64 &584.33 &2938.26 &1 &1 &1\\ 
      2  & 89.17 & 427.23 &2052.19 &1.36 &1.36 &1.43\\
      4 & 65.53  & 294.98 &1477.96 &1.86 &1.98 &1.99\\
      8 & 38.92 & 193.80 &853.67 &3.13 &3.02 &3.44 \\
      16 & 22.03 & 113.82 &511.41 &5.52 &5.13 & 5.74\\
      32 & 17.25 & 68.99 &308.04 &7.05 &8.47 & 9.54\\
      64 & 14.54 & 39.53 & 189.77 &8.37 &14.78 &15.48\\
      128 & 11.38 & 24.39 & 115.12 &10.69 &23.96 & 25.52 \\ \hline
    \end{tabular}
    \begin{tabular}[width = \textwidth]{c c c c c c c}
     \\
    \end{tabular}
    \quad \\
    \begin{tabular}[width = \textwidth]{c c c c c c c}
    \multicolumn{1}{c}{$h =1/256, \quad t =1$}\\ \hline
      $\tau$  & $1/256$ &  $1/512$ & $1/1024$ & $1/256$ &  $1/512$ & $1/1024$\\ \hline
      Number of threads  &  Wall time [s]  & & & Speedup ratio       \\
      1   & 850.10 &1565.33 &2938.26 &1 &1 &1 \\ 
      2  & 562.99 & 1067.04 &2052.19 &1.51 &1.47 &1.43 \\
      4 & 413.52  & 792.98 &1477.96  &2.06 &1.97 &1.99\\
      8 & 268.82 & 442.07 &853.67  &3.16 &3.54 &3.44\\
      16 & 163.40 & 257.47 &511.41 &5.20 &6.08 &5.74\\
      32 & 104.06 & 175.03 &308.04 &8.17 &8.94 &9.54 \\
      64 & 82.58 & 123.47 & 189.77 &10.29 &12.67 &15.48\\
      128 & 72.17 & 89.76 & 115.12 &11.78 &17.44 &25.52 \\ \hline
    \end{tabular}
    
\end{table}

\subsection{Barry \& Mercer's problem}
Anothrer benchmark problem for Biot's model is Barry \& Mercer's model. It describes the behavior of a rectangular domain with a sharp point source. For every boundary, liquid are drained and assume zero tangential displacement. The point source is of sine wave function:
\begin{equation}
    g(t) = -0.5 \delta_{(x_0,y_0)}sin(\beta t),
\end{equation}
where $\beta = \frac{(\lambda +2\mu)K}{ab}$ where the rectangular domain is $[0,a]\times[0,b]$, and $\delta_{(x_0,y_0)}$ is Dirac delta at point $(x_0,y_0)$. 

\begin{figure}
    \centering
    \includegraphics[width = 0.75\textwidth]{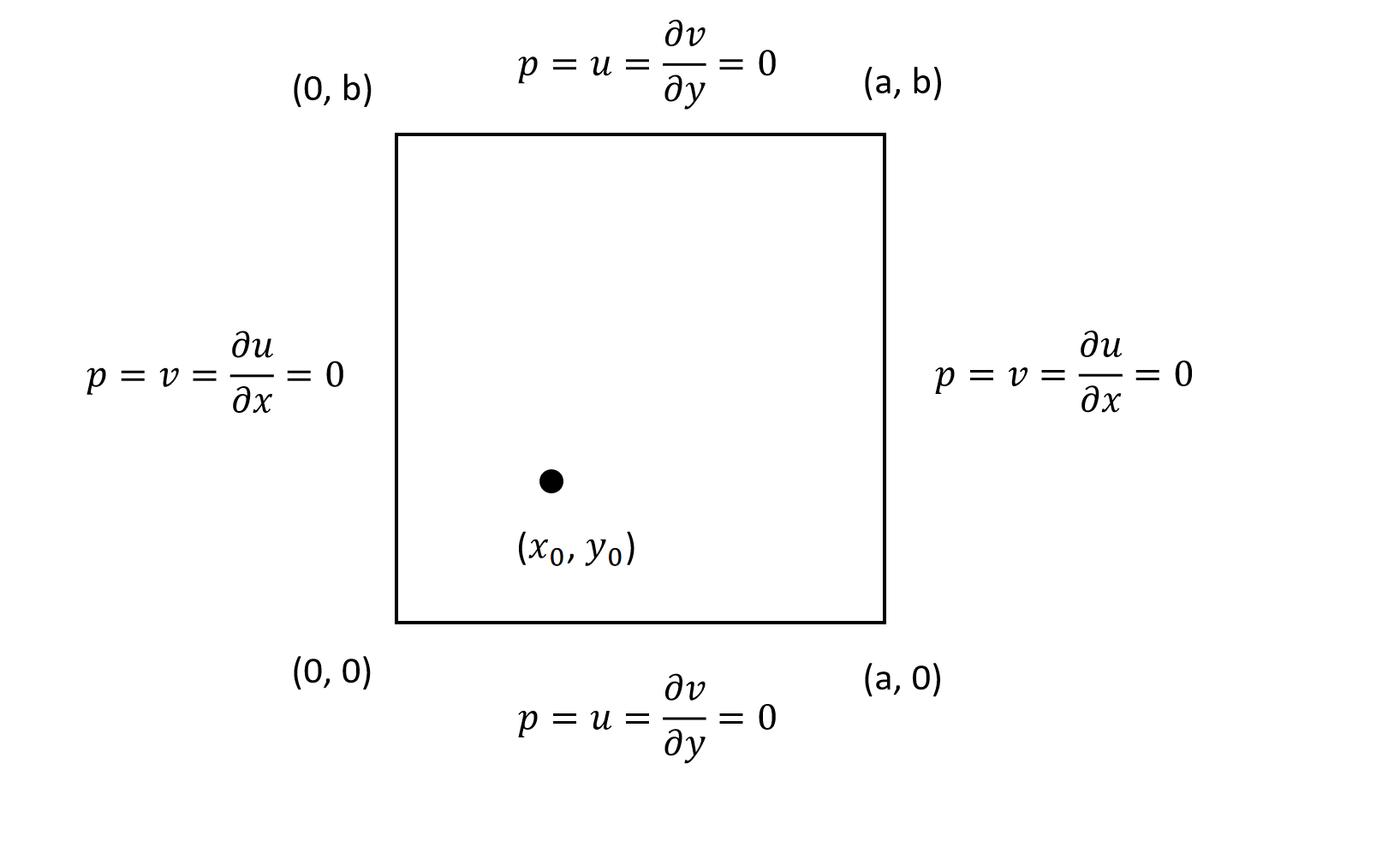}
    \caption{Domain and boundary conditions for Barry \& Mercer's Problem.}
    \label{fig:barry_domain}
\end{figure}

In this test case, choose $\Omega =[0,1]\times[0,1]$, and $x_0=y_0=1/4$. The material parameters are considered $E =10^5,\ \nu=0.1 \ \text{and}\ K=10^{-2}$. The boundary conditions are depicted in Fig. \ref{fig:barry_domain}.

\begin{figure}
\begin{subfigure}[t]{0.5\textwidth}
           \centering
           \includegraphics[width=1\textwidth]{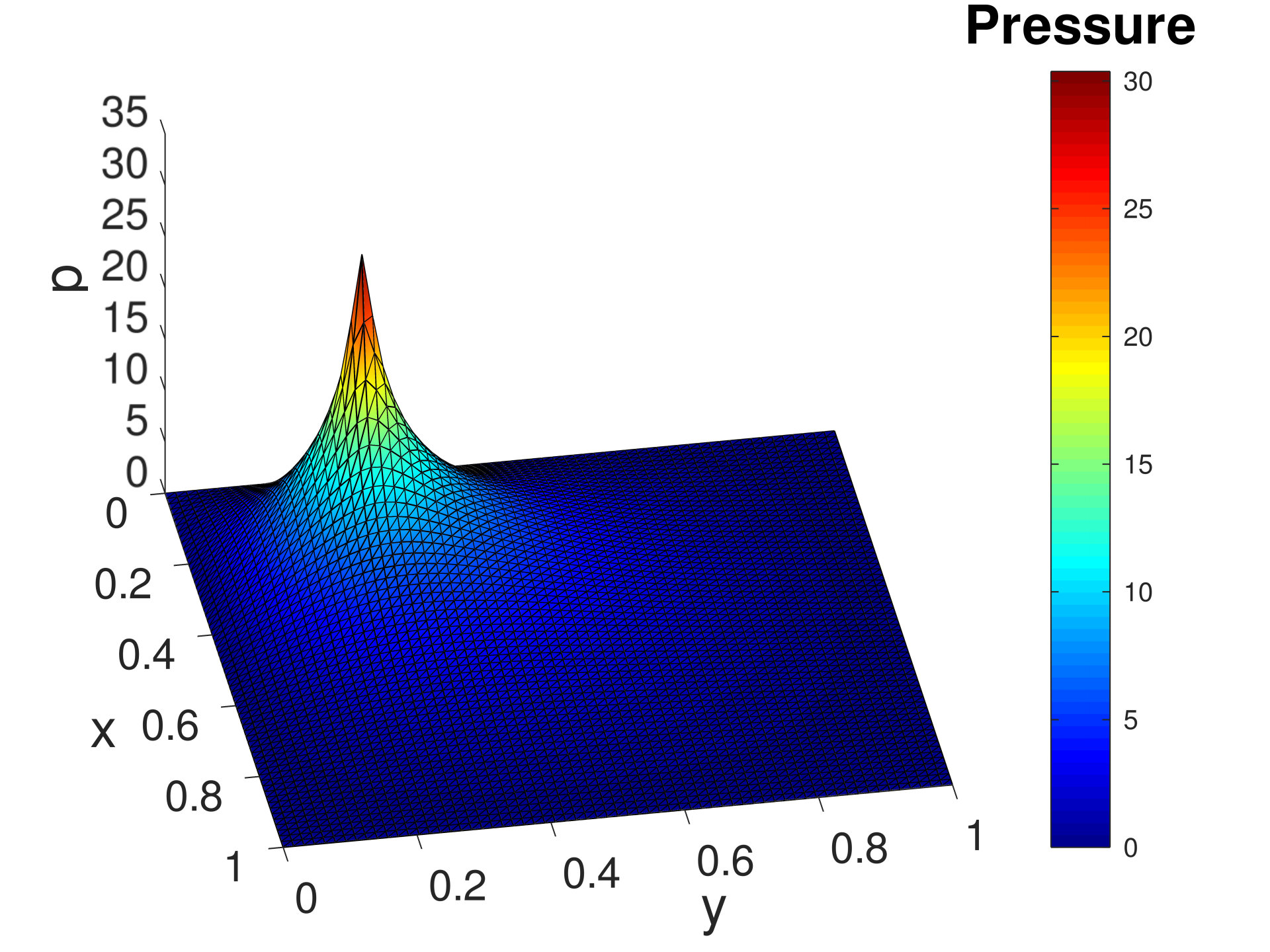}
            \caption{analytical solution}
            \label{fig:Case1}
    \end{subfigure}
    \begin{subfigure}[t]{0.5\textwidth}
           \centering
           \includegraphics[width=1\textwidth]{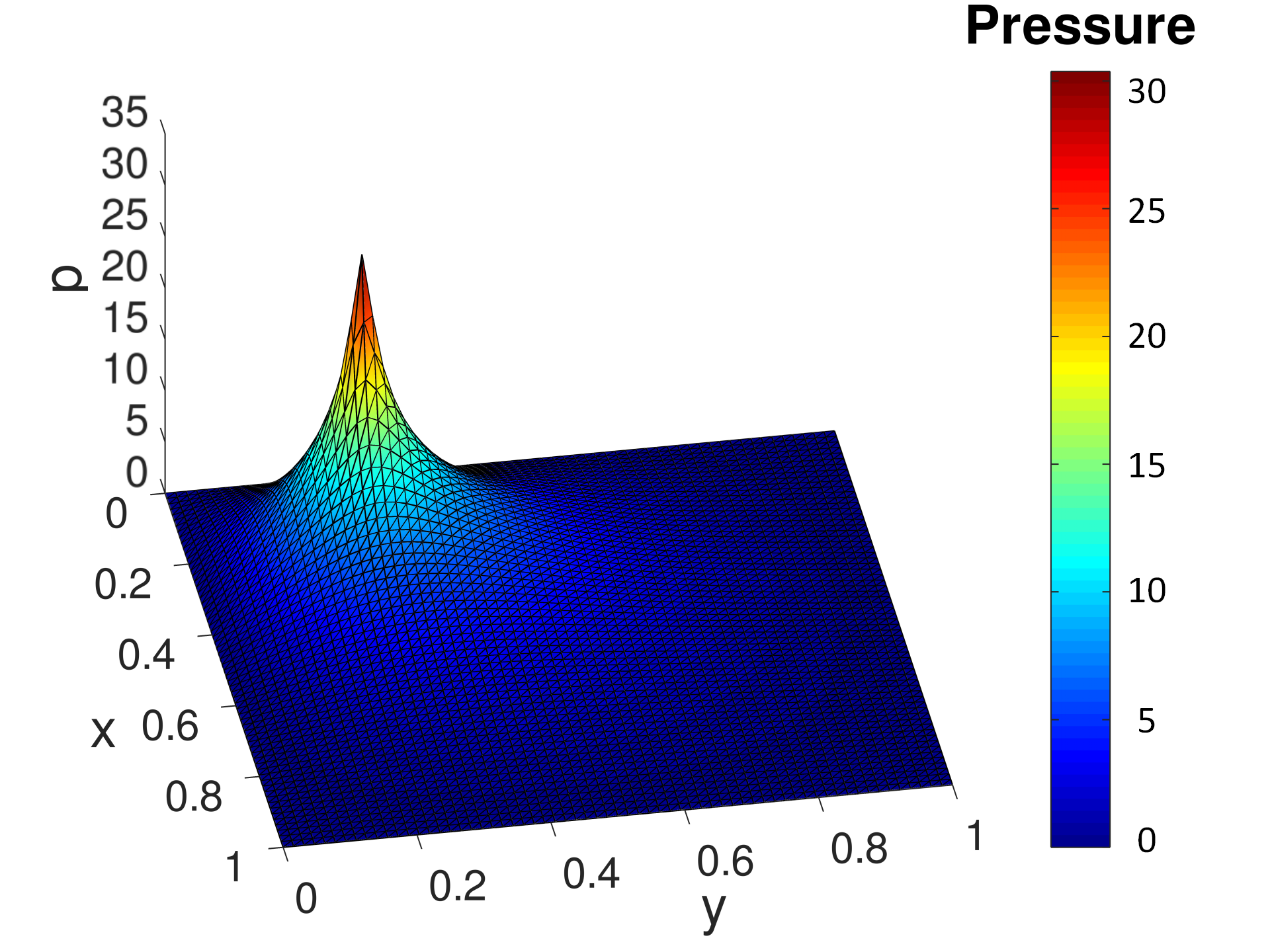}
            \caption{sequential backward Euler method}
            \label{fig:Case1}
    \end{subfigure}
    \begin{subfigure}[t]{0.5\textwidth}
            \centering
            \includegraphics[width=1\textwidth]{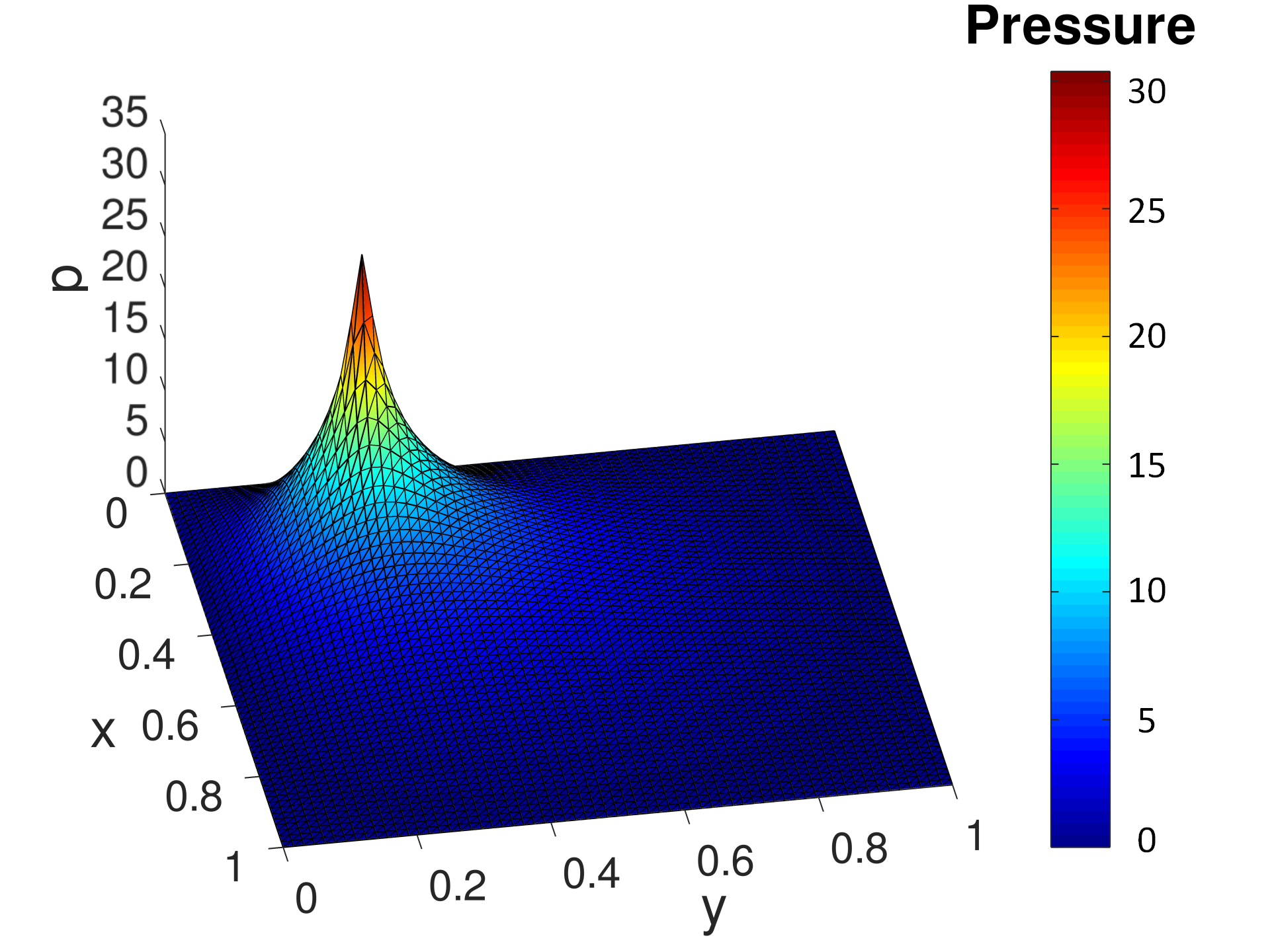}
            \caption{parallel 16 threads}
            \label{fig:b}
    \end{subfigure}
\begin{subfigure}[t]{0.5\textwidth}
           \centering
           \includegraphics[width=1\textwidth]{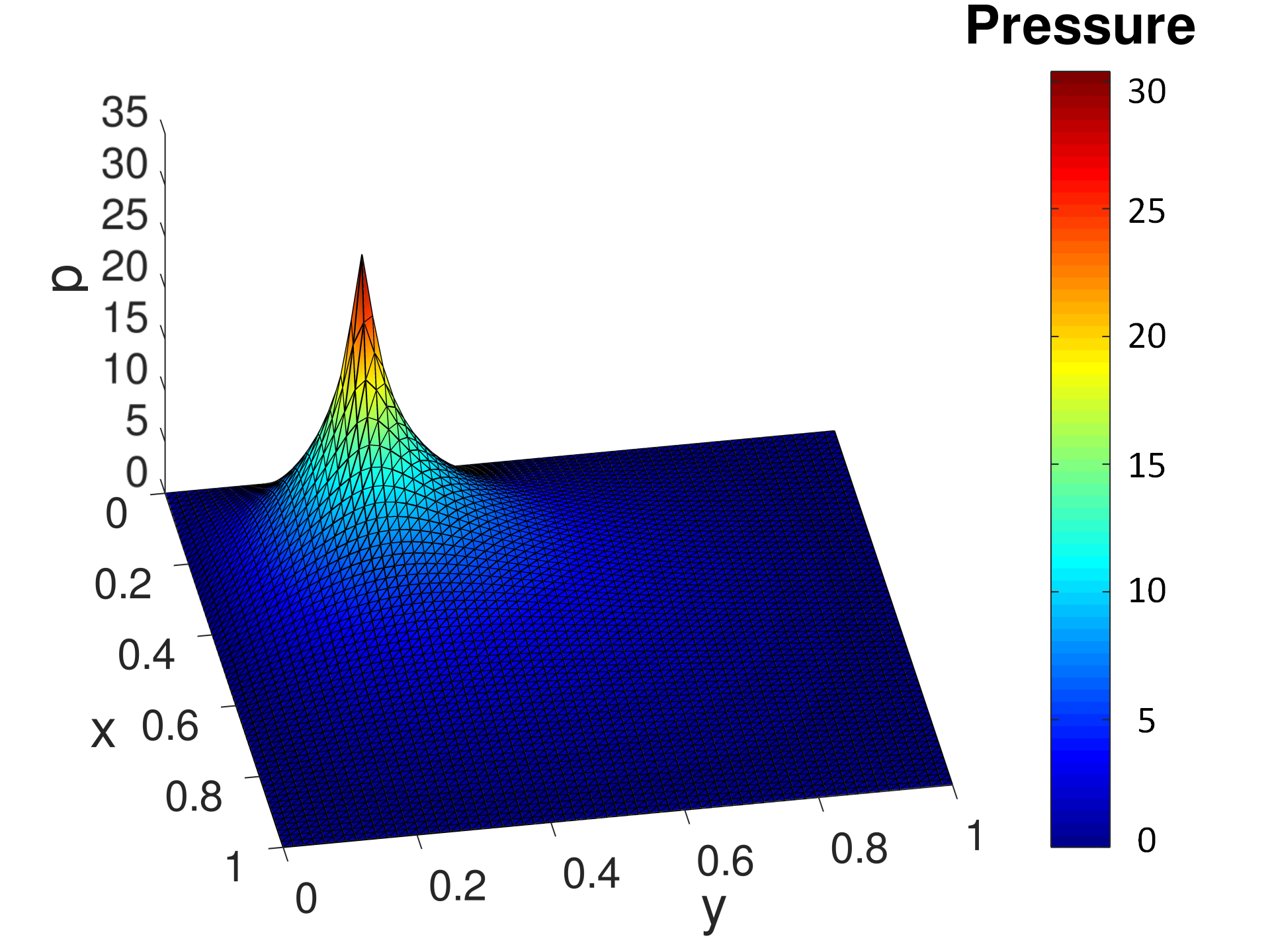}
            \caption{parallel 64 threads}
            \label{fig:Case1}
    \end{subfigure}
    
    \caption{Comparison of the solution of $p$ at $T=1$ with backward Euler method, parallel-in-time method and analytical solution.}
    \label{fig:T1_barry}

\end{figure}

Discrete the system in uniform right triangular grid with P1-P1 element and the spatial step length is $h_x=h_y=1/64$. The time step length $\tau =1/256T$ and the final time $T =1$. The analytical solutions is an infinite series and can be found in \cite{BarryMercy}. And the numerical solution in this report is presented in Fig. \ref{fig:T1_barry}. With stabilized P1-P1 scheme, the solution of p is still very smooth near the sharp point source. And the parallel-in-time method maintain this stability compared with sequential timestepping method. The Krylov subspace solver used in this case is the GMRES with preconditioner $\widetilde{P_2}$ in Eq. \eqref{p2}. The tolerance for the final residual in iterations is $10^{-10}$.

\begin{figure}

    \begin{subfigure}[t]{0.5\textwidth}
           \centering
           \includegraphics[width=\textwidth]{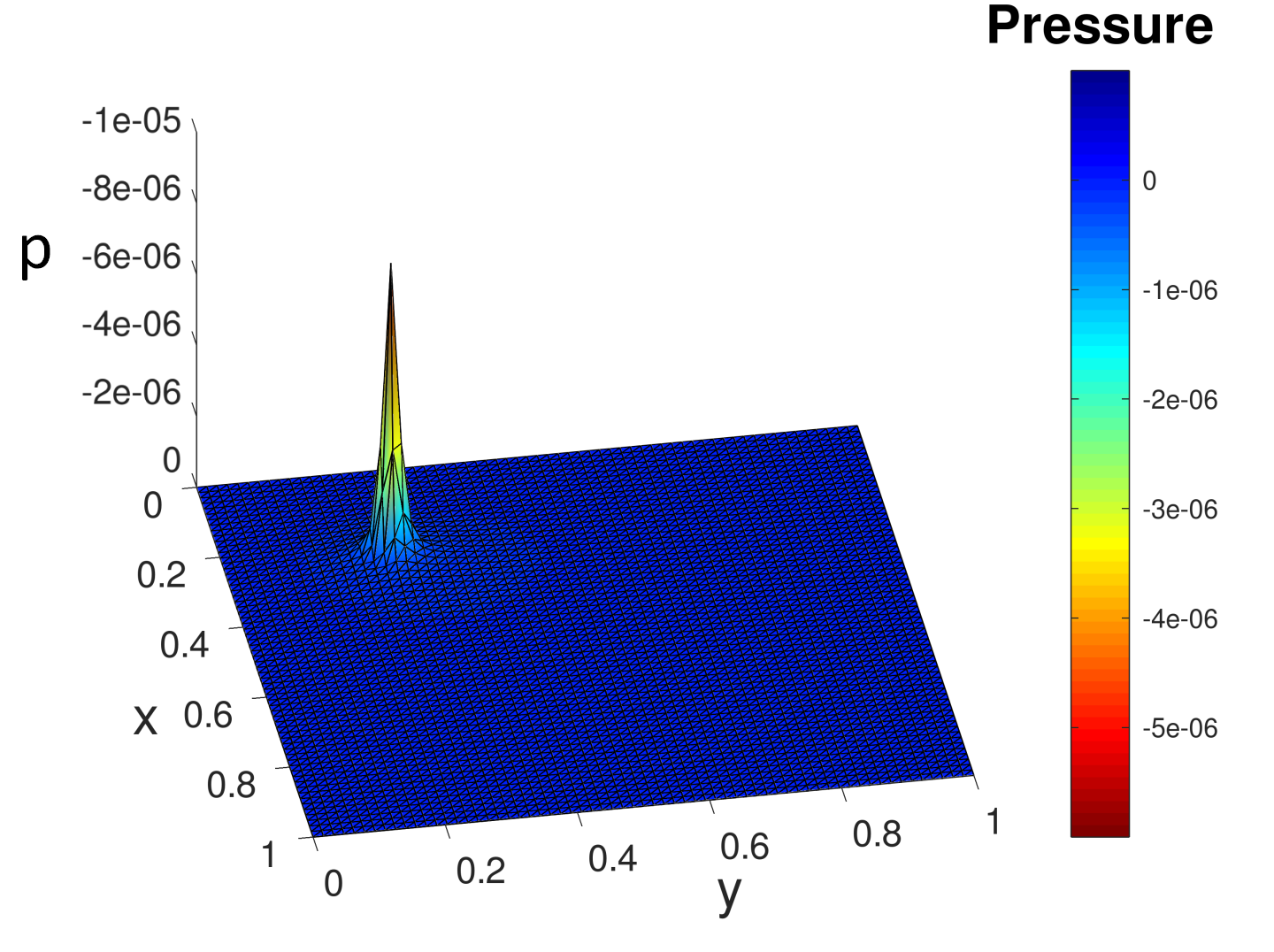}
            \caption{sequential (with stabilization term)}
            \label{fig:Case1}
    \end{subfigure}
    \begin{subfigure}[t]{0.5\textwidth}
           \centering
           \includegraphics[width=\textwidth]{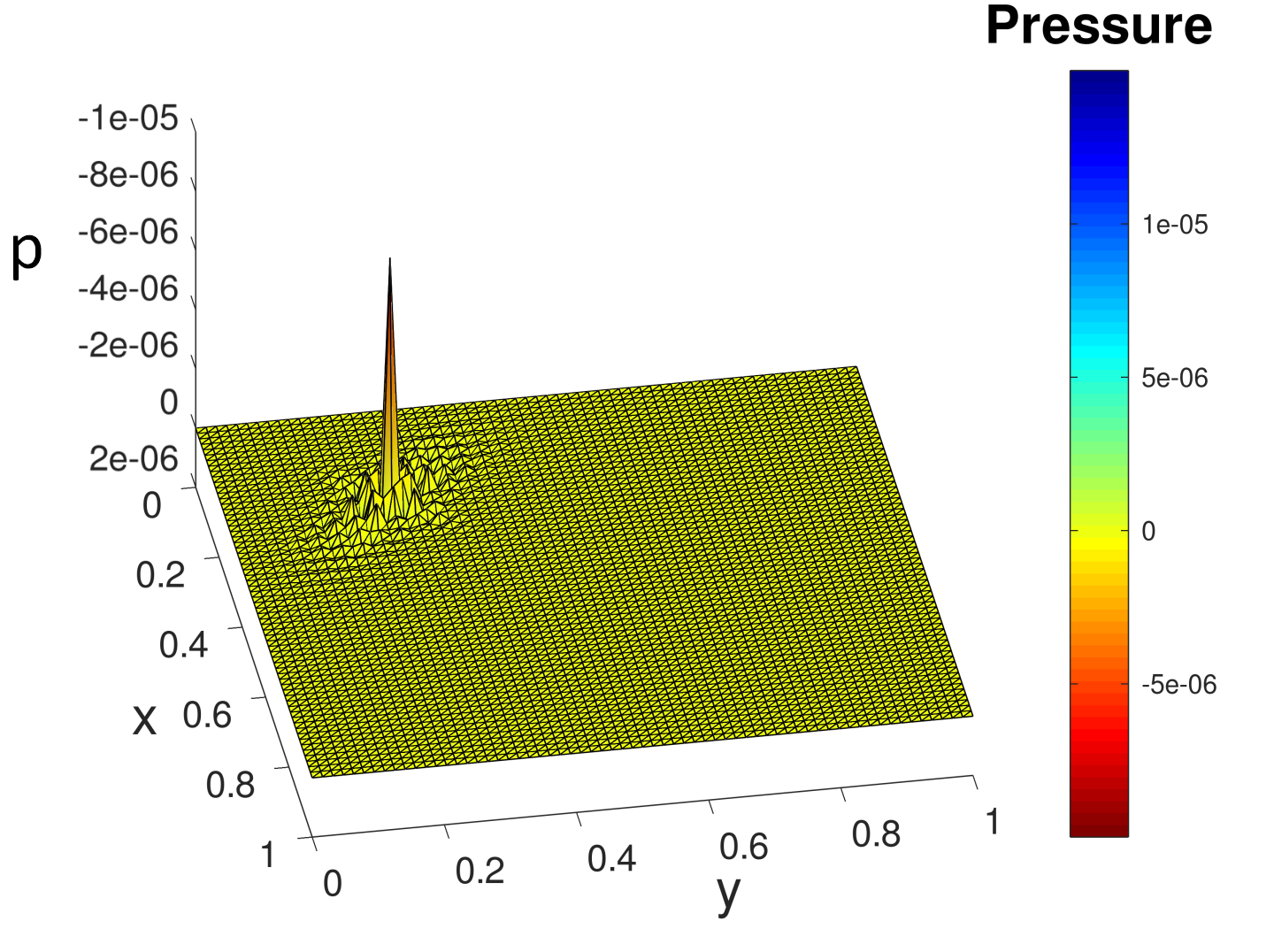}
            \caption{sequential (without stabilization term)}
            \label{fig:Case1}
    \end{subfigure}
    
    \begin{subfigure}[t]{0.5\textwidth}
            \centering
            \includegraphics[width=\textwidth]{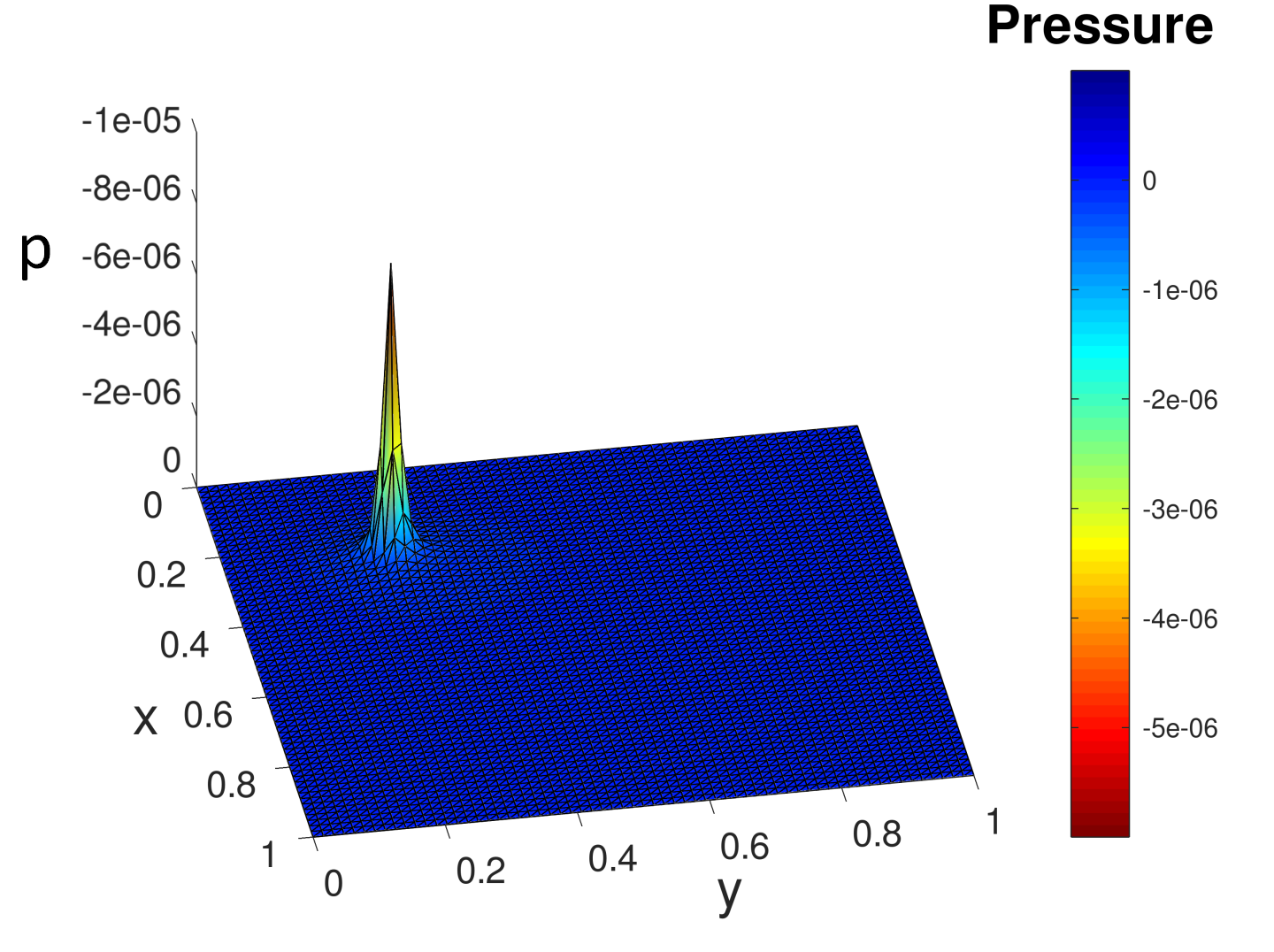}
            \caption{parallel 16 threads (with stabilization term)}
            \label{fig:b}
    \end{subfigure}
    \begin{subfigure}[t]{0.5\textwidth}
            \centering
            \includegraphics[width=\textwidth]{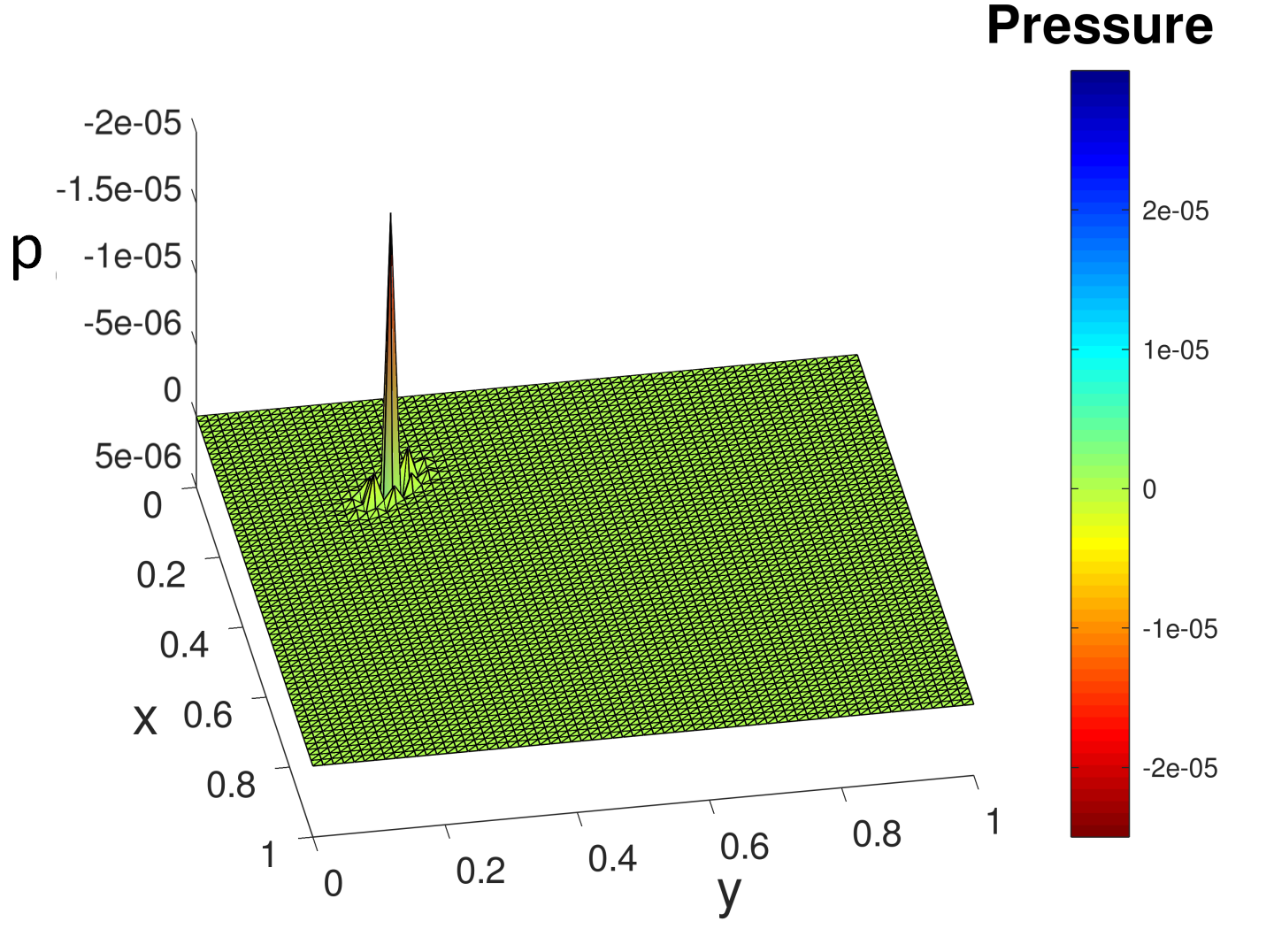}
            \caption{parallel 16 threads (without stabilization term)}
            \label{fig:b}
    \end{subfigure}
    
\begin{subfigure}[t]{0.5\textwidth}
           \centering
           \includegraphics[width=\textwidth]{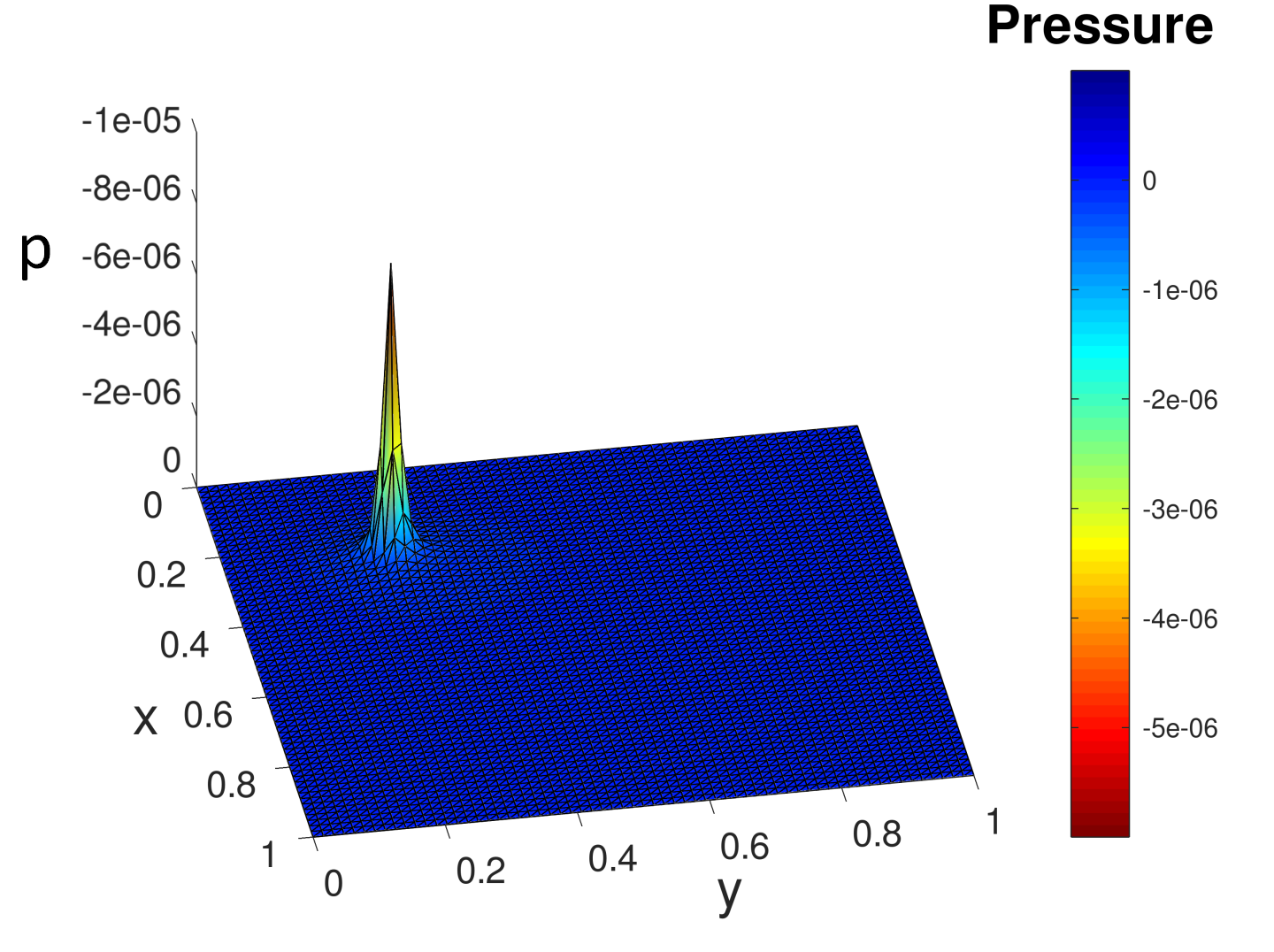}
            \caption{parallel 64 threads (with stabilization term)}
            \label{fig:Case1}
    \end{subfigure}
    \begin{subfigure}[t]{0.5\textwidth}
           \centering
           \includegraphics[width=\textwidth]{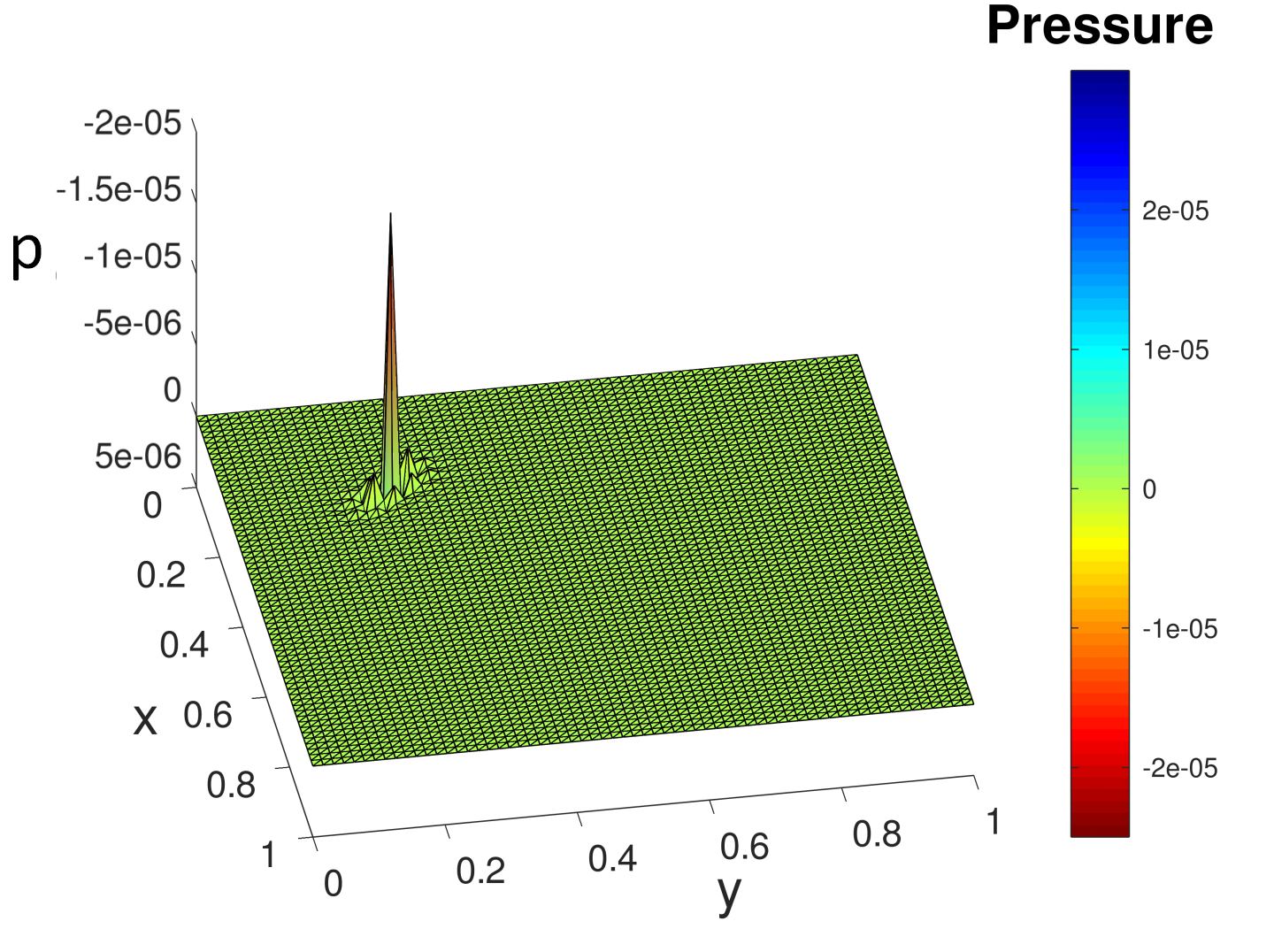}
            \caption{parallel 64 threads (without stabilization term)}
            \label{fig:Case1}
    \end{subfigure}
    
    \caption{Comparison of numerical solution and error of $p$ at $T=10^{-4}$ with backward Euler method and the parallel-in-time method.}
    \label{fig:barrystabornot}

\end{figure}

To further study the possible non-physical oscillation, choose $K=10^{-6}$, final time $T=10^{-4}$ and time step length $\tau =T/16$. The numerical solution computed by schemes with and without stabilization term are presented in Fig. \ref{fig:barrystabornot}. It is observed that, without stabilization term the parallel-in-time method outputs different results with the one given by sequential method, however, fortunately with the stabilized P1-P1 scheme, the solution given by the parallel-in-time method is still stable and consistent to the solution computed with sequential method.

\subsection{Mandel's Problem}
\label{sec:mandel}
Considering another bench problem for Biot's model - Mandel's Problem which is more close to the engineering application. At time $t=0$, an uniform vertical load of magnitude $2F$ is applied and an equal but upward force applied to the bottom plate. The magnitude of the force is constant during the simulation time. The domain is stress-free and free to move at boundary $x= \pm a$. Gravity and any other body force is neglected.

\begin{figure}

    \begin{subfigure}[t]{0.5\textwidth}
           \centering
           \includegraphics[width=\textwidth]{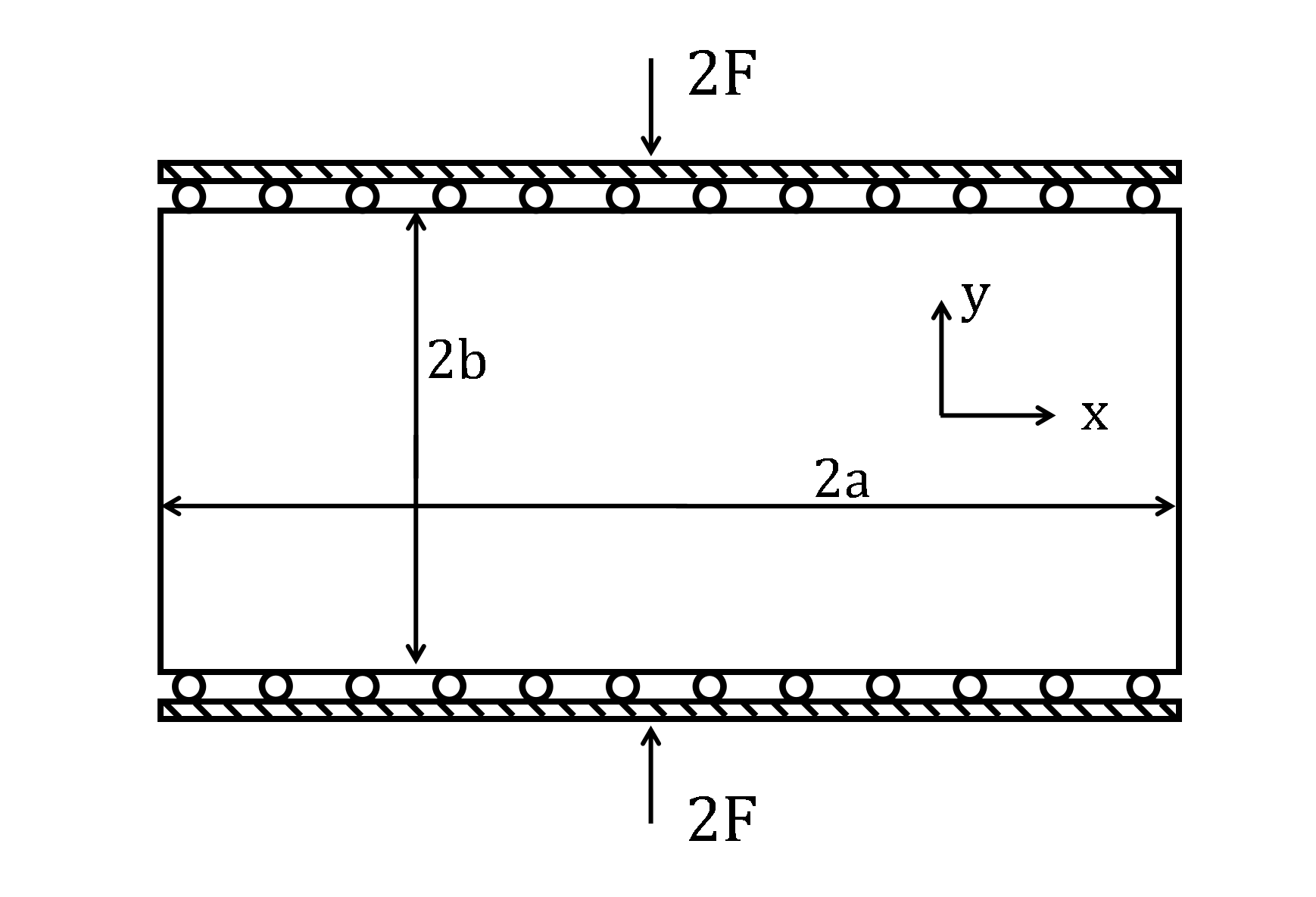}
            \caption{Physical domain for Mandel's problem }
            \label{fig:Case1}
    \end{subfigure}
    \begin{subfigure}[t]{0.5\textwidth}
            \centering
            \includegraphics[width=\textwidth]{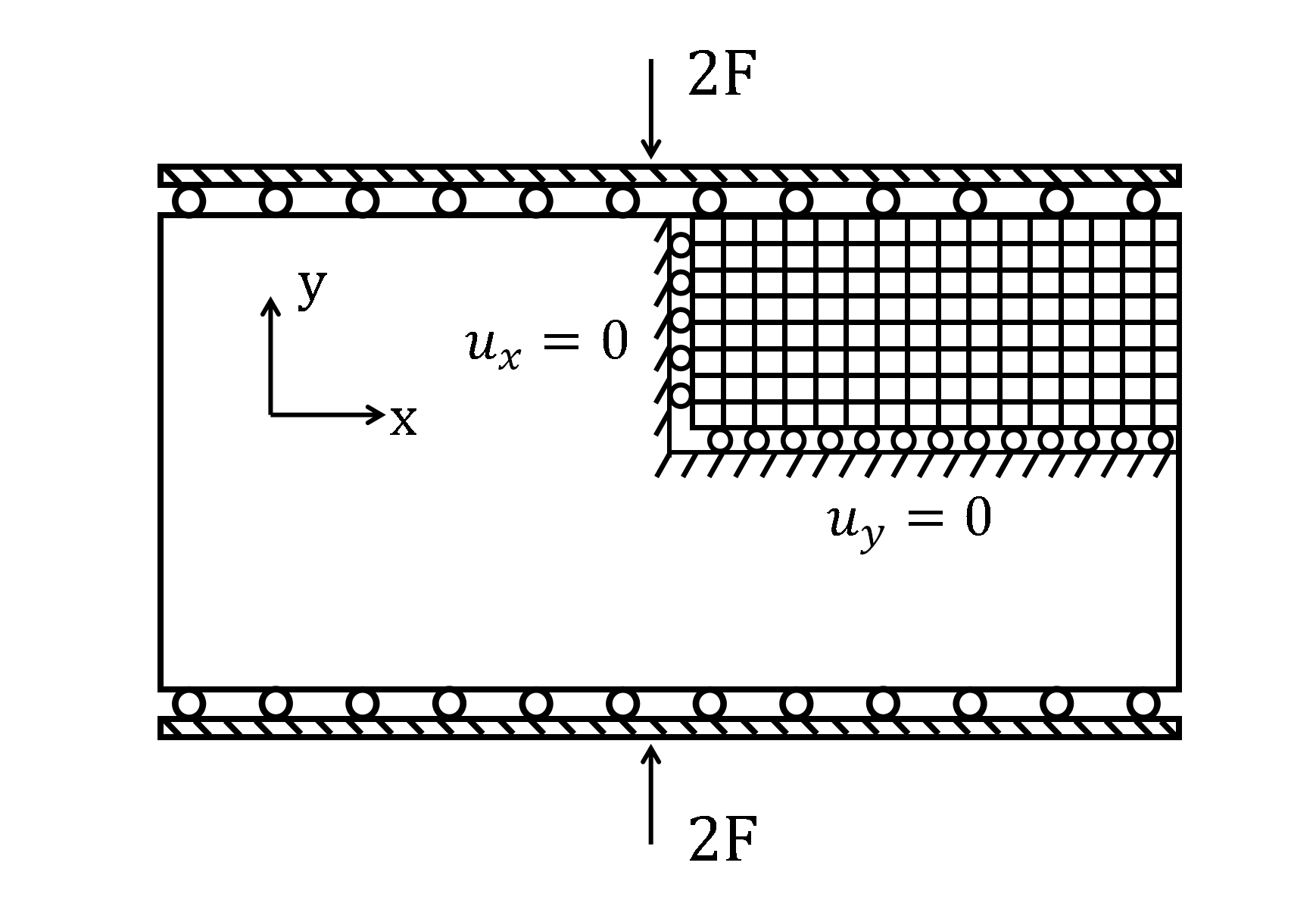}
            \caption{Quarter domain for Mandel's problem}
            \label{fig:b}
    \end{subfigure}

    \caption{Domain for Mandel's problem.}
    \label{fig:mandeldomain}
\end{figure}

Since the problem is symmetrical, it is resonable to use a quater domain and symmetical boundary conditions at inner face, see for details in Fig. \ref{fig:mandeldomain}. For the displacement of rigid plate on the top boundary, it can be enforced by constrained equations and the vertical displacement $u_y(b,t)$ is constant along x direction.
Due to the vertical load, an instantaneous and uniform pressure field will applied to the domain at initial time, and it is described analytically \cite{Mandelrevisit} and is used as initial condition:

\begin{equation}
\begin{aligned}
    p_0(x,y,0)&=\frac{FB(1+\nu_u)}{3a},
    \\
    \boldsymbol{u}(x,y,0)&=(\frac{F\nu_u x}{2G},\quad\frac{-FB(1-\nu_u) \nu)}{2Ga})^T,
    \end{aligned}
\end{equation}
where B is the Skempton coefficient and $\nu_u = \frac{3\nu+B(1-2\nu)}{3-B(1-2\nu)}$.
The boundary condition is described in Eq. \eqref{mandel boundary}
\begin{equation}
    \label{mandel boundary}
 \left\{
    \begin{aligned}
         p = 0,&\ \boldsymbol{\sigma} \cdot \boldsymbol{n}=0 \quad &\text{at}\ x=a \\
               \nabla p \cdot \boldsymbol{n} =0,&\ \boldsymbol{u}\cdot\boldsymbol{n} =0 \quad &\text{at}\ x=0 \\
                \nabla p \cdot \boldsymbol{n} =0,&\ \int_{\Gamma} \boldsymbol{\sigma} \cdot \boldsymbol{n}ds = -F\quad &\text{at}\ y=b
                \\
               \nabla p \cdot \boldsymbol{n} =0,&\
               \boldsymbol{u}\cdot\boldsymbol{n} =0 \quad &\text{at}\ y=0
       \end{aligned}    
    \right. .
\end{equation}
In this case, domain size is chosen as $a=1$ and $b=1$ . The material parameters $K=10^{-6}$, $B = 1$, $E=10^4 $ and $\nu = 0$, and therefore $\nu_u =0.5$. And the force on top boundary has magnitude of $F=1$. Using GMRES with preconditioner $\widetilde{P_1}$ in \eqref{p1} to solve the linear system in this case. The tolerance for the final residual in iterations is $10^{-7}$. 

\begin{figure}
    \centering
    \includegraphics[width=0.7\textwidth]{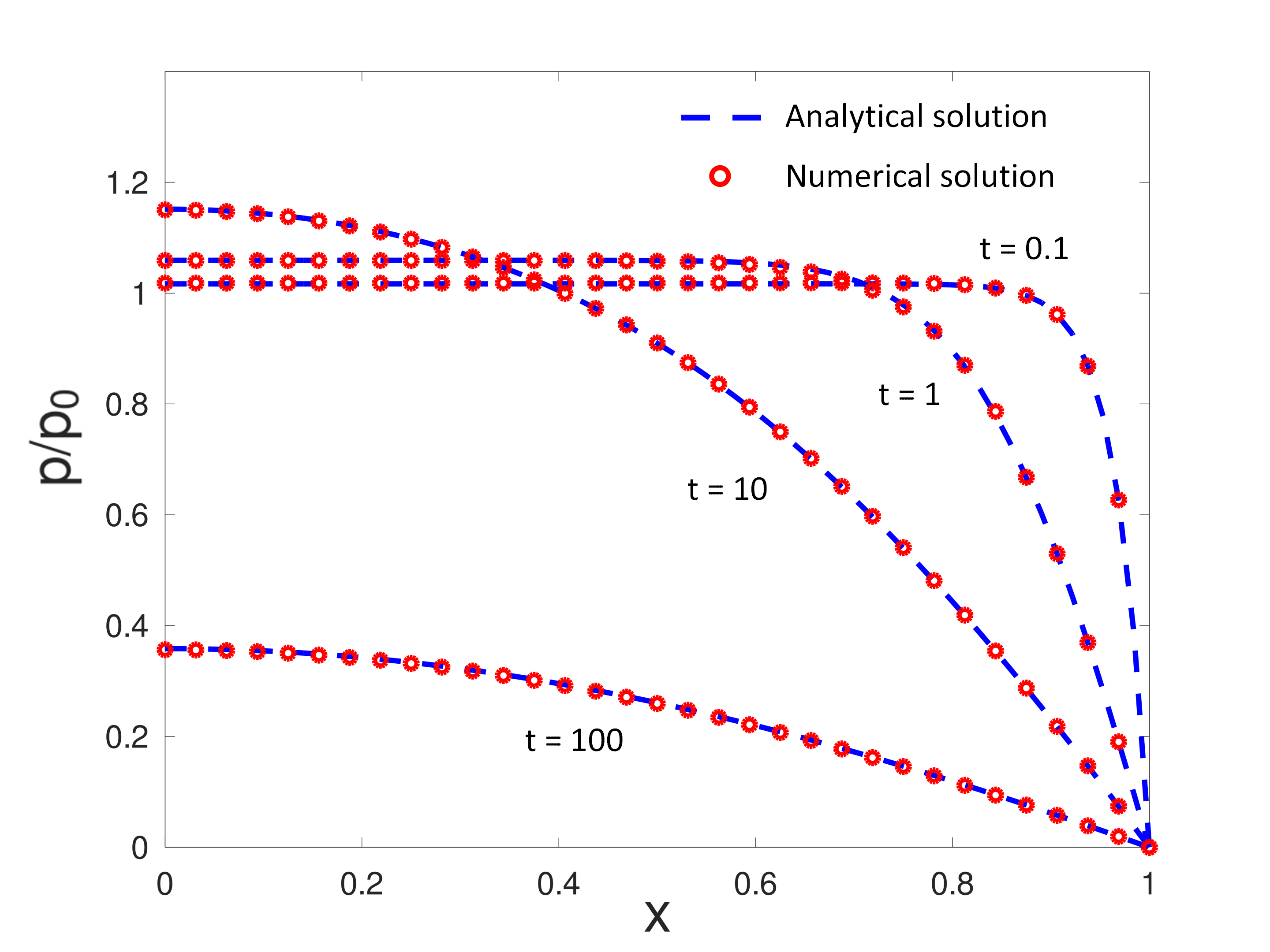}
    \caption{Comparison between numerical solution and analytical solution of Mandel's Problem, computed with parallel-in-time method by 16 threads.}
    \label{fig:Mandelnumerical}
\end{figure}

A comparison between the analytical solution and the numerical solution computed with paralle-in-time method by 16 threads is presented in Fig. \ref{fig:Mandelnumerical} when set the grid size $h_x=h_y=1/32$, and time step length $\tau=0.025$.

The numerical solution is compatible with the exact physical phenomenon, from the figure, it is observed that the pressure firstly increase and then decrease which is described as  Mandel–Cryer effect.

In this section, we are going to study the influence of different preconditioners and parallel architectures on the speedup performance. Set grid size $h_x=h_y=1/128$, time step length $\tau =0.1$ and total time length $T =100$.

\begin{figure}
    \centering
    \includegraphics[width=0.65\textwidth]{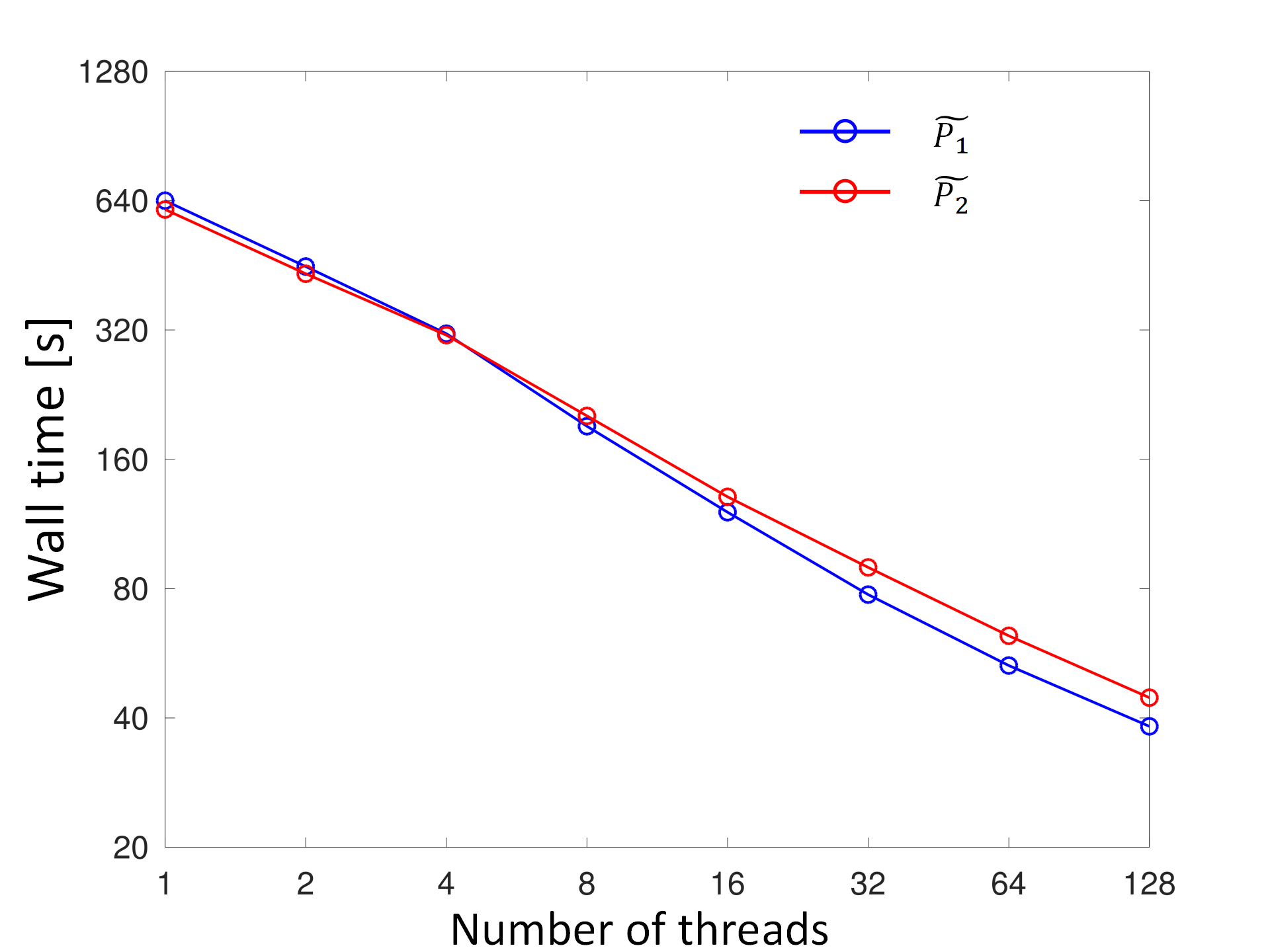}
    \caption{Wall clock time of parallel-in-time method with different preconditioners.}
    \label{fig:preconditionspeedup}
\end{figure}

To see the effect of different preconditioners on speedup rate, the wall clock times are presented in Fig. \ref{fig:preconditionspeedup}. In this case, the program runs in hybrid MPI-OpenMP architecture, every process conducts 4 threads. The results show that with both preconditioner the parallel-in-time method can achieve linear speedup performance. And the speedup performance with $\widetilde{P_1}$ is slightly better than that with $\widetilde{P_2}$. In addition, comparing the speedup ratio for Mandel's problem and that for cases in Section \ref{sec:trigono}, we can find that the practical speedup ratio varies between different cases. One possible reason is that different boundary conditions result in different matrix structure and influences the convergence process of Krylov subspace iterations.

\begin{table}
    \centering
    \caption{Wall time and speedup ratio of the parallel-in-time method with different preconditioners}
    \label{manprewalltimetable}
    \begin{tabular}[width = \textwidth]{c c c c c }
    \multicolumn{1}{c}{$\tau =0.1, \quad h_x =1/128,\quad h_y =1/128$}\\ \hline
      Preconditioner  & $\widetilde{P_1}$ &  $\widetilde{P_2}$ &$\widetilde{P_1}$ & $\widetilde{P_2}$  \\ \hline
      Number of threads  &  Wall time [s]  &\quad &Speedup ratio & \quad       \\
      1   & 640.20 &610.54 & 1 &1 \\ 
      2  & 449.78 & 432.71 & 1.43 &1.41\\
      4 & 313.83  & 311.37 & 2.04 & 1.96\\
      8 & 191.09 & 201.99 & 3.35 & 3.02\\
      16 & 120.52 & 130.89 & 5.33 & 4.66\\
      32 & 77.50 & 89.67  & 8.26 & 6.81\\
      64 & 52.98 & 62.15 & 12.08 & 9.82\\
      128 & 38.25& 44.61 & 16.74 & 13.69\\ \hline
    \end{tabular}
\end{table}

\begin{figure}
    \centering
    \includegraphics[width=0.7\textwidth]{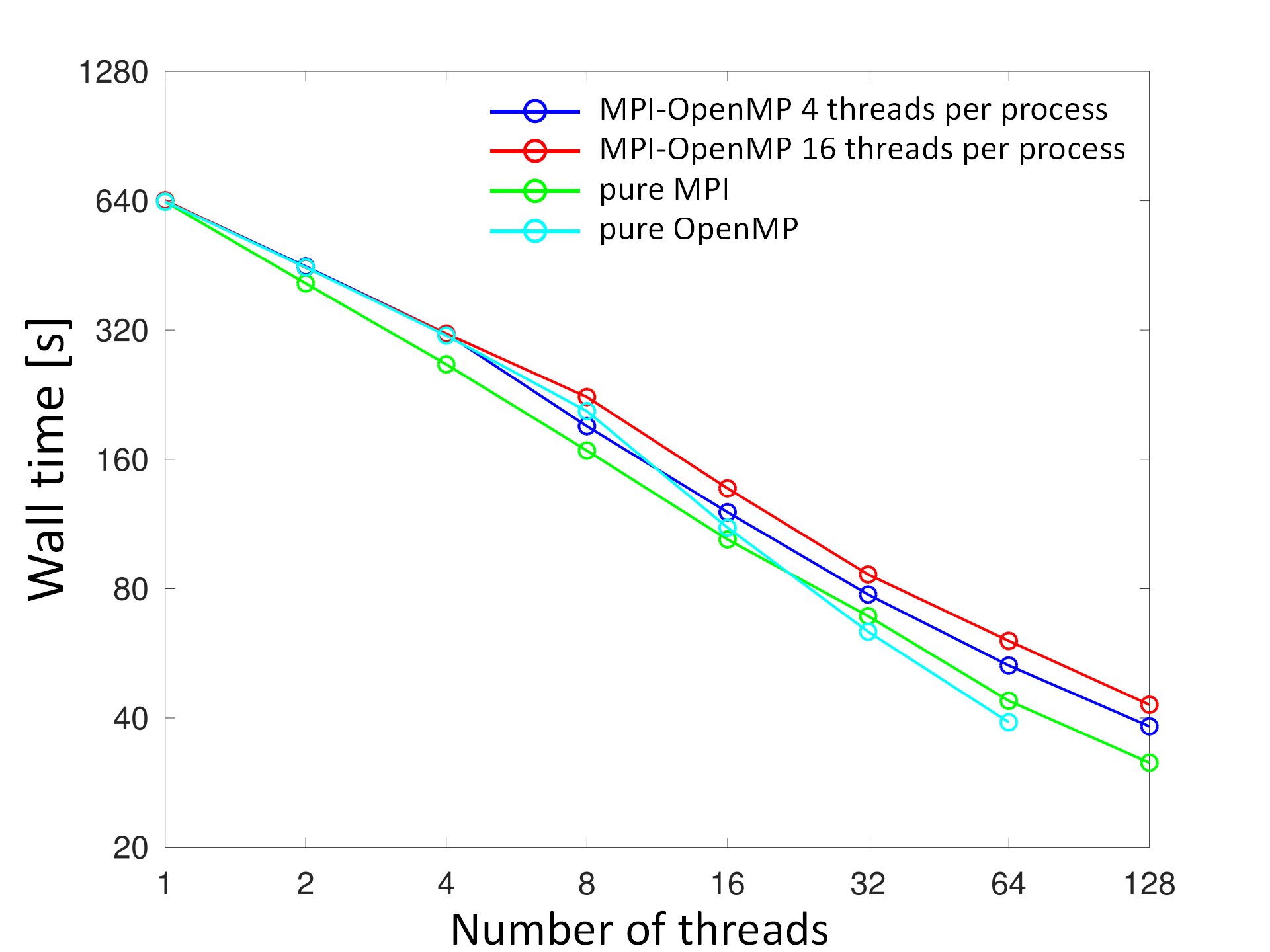}
    \caption{Wall clock time of parallel-in-time method with different architecture when $h_x=1/128,\ h_y=1/128,\ \tau =0.1$.}
    \label{fig:mandelarchi}
\end{figure}

To study the effect of different parallel architectures, run the program in pure MPI architecture, MPI-OpenMP hybrid architecture, and pure OpenMP architecture, wall-clock times are presented in Fig. \ref{fig:mandelarchi}. The preconditioner used here is $\widetilde{P^1}$. In the HPC platform where the test cases run, one computing node has two 32-core CPUs, so the pure shared memory parallel architecture with OpenMP has an upper limit for the number of threads of 64. From the results, it can be observed that the linear speedup performance is independent of the parallel architecture; the pure OpenMP architecture has the highest speedup ratio as the number of threads increases; however, when the number of threads increases from 1 to 8, the computing will slightly slow down. When applying pure MPI architecture, each CPU runs only one process, and one process only conducts one thread. Since point-to-point communication with MPI is less efficient than shared memory architecture in one CPU, the wall-clock time of pure MPI communication between processes in one CPU is meaningless to present here. With pure MPI architecture, there is no slow-down effect when the number of threads increases from 1 to 8, which results in a shorter wall-clock time. The speedup ratio in the later stage is similar to the hybrid MPI-OpenMP architecture. However, in practice, using only one core per CPU is very expensive while leaving the other cores resting. This test case proves that the algorithm can be adapted flexibly according to the hardware architecture of different HPC platforms and can consistently achieve a satisfying speedup performance.
\begin{table}
    \centering
    \caption{Wall time and speedup ratio of the parallel-in-time method with different parallel architecture}
    \label{manarcwalltimetable}
    \begin{tabular}[width = \textwidth]{c c c c c}
    \multicolumn{1}{c}{$\tau =0.1, \quad h_x =1/128,\quad h_y =1/128$}\\ \hline
      Parallel Architecture  & MPI-OpenMP  &  MPI-OpenMP  & pure MPI & pure OpenMP  \\
      \quad & 4 threads/ process & 16 threads/ process & \quad & \quad \\ \hline
      Number of threads  &  Wall time [s]          \\
      1   & 640.20&	641.16&	636.68&	637.15 \\ 
      2  & 449.78&  447.20 &	410.96 & 447.86 \\
      4 & 313.83 & 313.91 &	266.20&	310.69 \\
      8 & 191.09&	223.43&	167.69&	207.11 \\
      16 & 120.52 & 136.96 &104.15 &110.97\\
      32 & 77.50&	86.33&	69.06&	63.50\\
      64 & 52.98&	60.46&	43.85&	39.11\\
      128 & 38.25&	42.95&	31.48 & Not Available \\ 
      \quad \\
      Number of threads  &  Speedup ratio \\
      1   & 1 &1 & 1 &1 \\ 
      2  & 1.43 & 1.43 & 1.55 & 1.42 \\
      4 & 2.04  & 2.04  & 2.39 & 2.05\\
      8 & 3.35 & 2.87  & 3.78 & 3.08 \\
      16 & 5.33 & 4.68 &6.11 &5.74\\
      32 & 8.26& 7.43 &9.22 &10.03 \\
      64 & 12.08 & 10.61 &14.52 &16.29\\
      128 & 16.74 & 14.91 & 20.22 & Not Available \\ 
      \hline
    \end{tabular}

\end{table}

\clearpage

\section{Conclusion}\label{sec:conclusion}
In this paper, we propose the parallel-in-time method based on preconditioner for Biot's problem, and conduct numerical experiments to demonstrate that the parallel-in-time method preserves the consistency and accuracy of the original sequential method. Moreover, we show that the parallel-in-time method achieves linear speedup on HPC platforms, and maintains this performance as the number of threads increases from 1 to 128, which indicates its potential for massively parallel computing. In the future, it is an important topic to analytically investigate the consistency of the solutions of time-parallel methods with traditional serial algorithms. Moreover, the parallel-in-time method does not conflict with spatial parallel methods such as domain decomposition, and  developing hybrid time-space parallel methods to achieve higher parallel efficiency would be one of the promising research directions.

\clearpage
\backmatter

\bmhead{Acknowledgments}
We thank National Supercomputing Center in Shenzhen (Shenzhen, China) for providing the internship opportunity and computing resource required by this paper.

\clearpage
\bibliography{sn-bibliography}


\end{document}